\documentclass[12pt]{article}

\usepackage[utf8]{inputenc}
\usepackage{amsmath,amsfonts}
\usepackage{bbm}
\usepackage{amsthm}
\usepackage{setspace}
\usepackage{subcaption}
\usepackage{fancyhdr}
\usepackage{graphicx}
\usepackage{appendix}
\usepackage{authblk}
\usepackage{mathtools}
\usepackage{breqn}
\usepackage{mathrsfs}
\usepackage[round]{natbib}
\usepackage{xcolor}
\usepackage{url}
\usepackage{marvosym}

\usepackage[left=2.5cm,right=2.5cm,top=2.5cm,bottom=2.5cm]{geometry}

\newtheorem{Theo}{Theorem}
\newtheorem{Lemma}{Lemma}

\newtheorem{Conj}{Conjecture}

\newcommand{\Esp} {\mathbbm{E}}

\newcommand{\Var}{\mathsf{Var}}

\newcommand{\Ind}{\mathbbm{1}}
\newcommand{\Pro}{\mathbbm{P}}

\newcommand{\Polish}{(\mathcal{X},\mathcal{B}_{\mathcal{X}})}
\newcommand{\UnitInt}{([0,1],\mathcal{B}_{[0,1]})}
\newcommand{\SeqRV}{X^{(\infty)}}
\newcommand{\lengthVar}{\{v_i\}_{i=1}^{\infty}}
\newcommand{\locations}{\{z_i\}_{i=1}^{\infty}}
\newcommand{\jumps}{\{w_i\}_{i=1}^{\infty} }

\def\simiid{\stackrel{\mbox{\scriptsize{iid}}}{\sim}}

\providecommand{\keywords}[1]
{
  \small	
  \textit{Key words:} #1

}

\title{On a divergence-based prior analysis of stick-breaking processes}
\author{José A. Perusquía$^1$, Mario Diaz$^2$\Cross, Ramsés H. Mena$^2$}
\affil{\small{Department of Mathematics, Faculty of Sciences, UNAM$^1$\\
Department of Probability and Statistics, IIMAS, UNAM$^2$}}
\date{}

\begin{document}
\setstretch{1.2}

\maketitle
\begin{abstract}
The nonparametric view of Bayesian inference has transformed statistics and many of its applications. The canonical Dirichlet process and other more general families of nonparametric priors have served as a gateway to solve frontier uncertainty quantification problems of large, or infinite, nature. This success has been greatly due to available constructions and representations of such distributions. Hence, understanding their distributional features and how different random probability measures compare among themselves is a key ingredient for their proper application. In this paper, we analyse the discrepancy among some relevant nonparametric priors. Initially, we compute the mean and variance of the random Kullback-Leibler divergence between the Dirichlet process and the geometric process. Subsequently, we extend our analysis to encompass a broader class of exchangeable stick-breaking processes, which includes the Dirichlet and geometric processes as extreme cases. Our results establish quantitative conditions where all the aforementioned priors are close in total variation distance. 
\end{abstract}


\section{Introduction}
\label{introduction}
Bayesian nonparametric models have become increasingly popular within many fields of research, due to their flexibility and the computational advancements that have made these models no longer intractable. However, the mathematical and computational complexity of specifying prior distributions on infinite-dimensional spaces, clearly posses a challenging task \cite[\textit{see e.g.}][]{hjort_holmes_muller_walker_2010,ghosal_van_der_vaart_2017}.  
With this in mind,  a preferred approach by practitioners, to define Bayesian nonparametric priors, is the one based on stick-breaking processes. Indeed, though not often used in its full generality, by modulating the distributional and dependence assumptions of the corresponding length variables (sometimes also named stick or break variables), the Dirichlet process \citep{ferguson,Sethuraman94} as well as the general class of normalised homogeneous completely random measures \citep{FLNNPT16}, have stick-breaking representations. That said, due to their mathematical and computational complexity, one could say that there is a sub-usage of dependent length variable stick-breaking representations. Recently, a general class of stick-breaking processes with exchangeable length variables was introduced by \cite{Leyva_2021}, having as particular cases the independent (\textit{e.g.} Dirichlet process) and fully dependent (\textit{e.g.} geometric process) length variable cases, which is a novel, robust and general strategy to build tractable nonparametric priors. Thus, it is of central interest understanding the differences among priors in this broader class.

Of course, more general, more robust and better suited stick-breaking processes could be used depending on the data and the application at hand. However, it is quite possible that these priors will lead to models which are both more mathematically and more computationally challenging. There is a clear trade-off between mathematical/computational ease and the usage of more general models, which needs to be carefully assessed when specifying the prior. In this work, we propose a divergence-based approach to compare stick-breaking priors. Specifically, we employ the Kullback-Leibler divergence to assess the impact of dependence in the length variables within the class of exchangeable stick-breaking priors introduced in \cite{Leyva_2021}. To the best of our knowledge, this kind of divergence-based analysis on stick-breaking priors has never done before. The closest approach to our methodology is perhaps \cite{watson_nieto_holmes} where the authors characterised the dispersion of some random probability measures using the Kullback-Leibler divergence. In their work, the main objective is to understand how far are the sample draws of such measures with respect to their baseline distribution, hence, allowing the authors to provide interesting insight on how to choose the truncation level and the parameters.

Our methodology allows us to unveil quantitative conditions under which the aforementioned priors exhibit close total variation distances. So that, we are able to share insight on the complexity cost of considering more complicated stick-breaking processes as well as to calibrate the gain when switch from one model to another. The remainder of the paper is organised as follows: In Section~\ref{preliminaries}, we present some preliminaries of stick-breaking processes, total variation distance, and Kullback-Leibler divergence. We compute the mean and variance of the Kullback-Leibler divergence between the Dirichlet and geometric processes in Section~\ref{KL_DIR_GEOM}. In Section~\ref{exch_sbp}, we extend our analysis to encompass the class of Dirichlet-driven exchangeable length-variables and study their Kullback-Leibler divergence with respect the geometric process. We discuss some final remarks and future lines of research in Section~\ref{conclusions}.

\section{Preliminaries}
\label{preliminaries}

\subsection{Stick-breaking processes}

Bayesian statistics can be closely linked and motivated by the mathematical notion of exchangeability \citep{schervish1995theory}. Within the Bayesian nonparametric realm, the representation theorem for exchangeable random variables \cite[][]{deFinetti31,Hewitt1955} establishes that an infinite sequence of random variables, henceforth denoted as $\SeqRV:=(X_i)_{i\geq 1}$, defined over a Polish space endowed with its Borel $\sigma$-algebra, $\Polish$,  is exchangeable if and only if it exists a measure $Q$ on the space of all probability measures on $\mathcal{X}$, denoted by $\mathcal{P}_{\mathcal{X}}$, such that for any $n\geq 1$
\begin{eqnarray*}
X_i\mid P &\simiid& P, \quad i=1,\ldots,n,\\
P &\sim& Q.\nonumber 
\end{eqnarray*}

In words, given the unknown shape of the distribution depicting the random variables $X_i$'s, they are independent and identically distributed (iid), provided that such form is random and uniquely modelled by $Q$.  Within this framework, $Q$ is unique, usually referred to as the de Finetti's measure, and plays the role of the prior distribution for $P$.

Clearly, specifying $P$ and its law $Q$, lies at the core of Bayesian nonparametric methods and there are at least two broad approaches to doing so. The first approach aims to directly specify $Q$, two of the most common ways to do it is through the
Kolmogorov's extension theorem for projective finite dimensional distributions or, indirectly through the predictive distribution that completely characterises the law of $\SeqRV$. As for the second approach, it aims to provide a
suitable stochastic representation of $P$. In this scenario two well-known approaches are based on a suitable transformation of a certain class of stochastic processes or based on the well-known representations via species sampling models. In this paper we centre our attention on the latter. Now, if $P$ is a species sampling process with base measure $P_0$,  then $P$ has the following a.s. representation
\begin{eqnarray*}
P&=&\sum_{n=1}^{\infty}w_n\delta_{z_n} + \left(1-\sum_{n=1}^{\infty}w_n\right)P_0,
\end{eqnarray*}
where $\sum_{n=1}^{\infty}w_n\leq 1$ almost surely (a.s.) and the weights, $\jumps$, are independent of the locations, $\locations$, which are a collection of iid random variables drawn from a non-atomic distribution $P_0=\Esp(P)$. 

In the case where $\sum_{n=1}^{\infty}w_n= 1$ a.s., we say that $P$ is proper and we can then consider the case of Bayesian nonparametric priors known as stick-breaking processes, whose weights are completely specified as
\begin{align}
\label{SBP}
    w_1=v_1, &&w_n=v_n\prod_{j=1}^{n-1}(1-v_j), \hspace{20pt}j>1,
\end{align}
where $\lengthVar$ $\in[0,1]$ are known as the length variables. Interesting particular cases are the well-known Dirichlet process when the $\lengthVar$ are independent and identically distributed as $Be(1,\theta)$ \citep{Sethuraman94} and the geometric process where $v_i\sim \delta_v$ with $v\sim Be(a,b)$, so that $w_n=v(1-v)^{n-1}$ \citep{Fuentes,Mena2011}. Here $Be(a,b)$ denotes the beta distribution with mean $a/(a+b)$. Of course, more general length variables, $\lengthVar$, can be chosen by letting $v_i\sim Be(a_i,b_i)$. In this case a sufficient condition to ensure a proper $P$ is that $\sum_{i=1}^{\infty}\Esp(\log(1-v_i))=-\infty$ or alternatively that $\sum_{i=1}^{\infty}\log(1+a_i/b_i)=\infty$ \citep{Ishwaran01}. 

Turning back our attention to the Dirichlet and the geometric process, it is compelling to remark that these two processes respectively represent examples of independent and fully dependent length variables. Indeed, they were studied apart until the recently developed theory in \cite{Leyva_2021}, where it was showed that if the length variables $\lengthVar$ are themselves exchangeable, and hence, driven by a random probability measure $\nu$ with base measure $\nu_0$, then, under certain conditions, we can recover both processes through the tuning of the tie probability $\rho_{\nu}=\Pro(v_1=v_2)=p(2)$. 

A practical way to define exchangeable length variables is to model $\nu$ via a nonparametric prior for which we are able to obtain a closed-form expression of the \textit{exchangeable partition probability function} (EPPF) $p(\cdot)$, which in turn  allows one to explicitly know the tie probability. Recall that for a given sequence of exchangeable random variables, governed by a discrete nonparametric prior, there is an associated distribution over the set of partitions induced by ties in the sequence \citep{pitman95}. In fact, for a given partition of the sequence, $(A_1,\ldots,A_k)$, the EPPF is a symmetric function that only depends on the cardinalities of the sets, {\it i.e.,} $(n_1,\ldots,n_k)$, with $n_j:=|A_j|$. As a direct example of this, if $\nu$ is a Dirichlet process with mass parameter $\beta$ and diffuse base measure $\nu_0=Be(1,\theta)$, then the tie probability is given by $\rho_{\nu}=\frac{1}{1+\beta}$ or equivalently, $\beta=\frac{1-\rho_{\nu}}{\rho_{\nu}}$. Hence, we can recover as limit cases the Dirichlet process when $\beta\rightarrow\infty$ (or $\rho_{\nu}\rightarrow0$) and the geometric process when $\beta\rightarrow0$ (or $\rho_{\nu}\rightarrow1$). 

\subsection{Total variation distance and Kullback-Leibler divergence}
\label{Subsection:TotalVariationDistance}

Let $P$ and $Q$ be two measures defined on a measurable space $(\varOmega,\mathscr{F})$. Recall that their total variation distance is defined as
\begin{equation*}
    \mathrm{TV}(P,Q) = \sup_{A\in\mathscr{F}}|P(A)-Q(A)|.
\end{equation*}
Total variation distance plays a key role in several statistical problems due to its properties, see, \textit{e.g.}, Chapter~7 of \cite{polyanskiy2022information}. For instance, total variation distance has a \textit{variational representation} that establishes that
\begin{equation*}
    \mathrm{TV}(P,Q) = \frac{1}{2} \sup_{f:\Omega\to[-1,1]} \big\lvert \Esp_{X \sim P}[f(X)] - \Esp_{X \sim Q}[f(X)] \big\rvert.
\end{equation*}
In particular, the previous equation implies that if $\mathrm{TV}(P,Q)$ is small, then any statistical inference of the form $\Esp[f(X)]$ with $f$ bounded would be similar for both $P$ and $Q$, {\it i.e.,}, $P$ and $Q$ are statistically indistinguishable.

Despite its usefulness, total variation distance is often challenging to compute in practice, making it necessary to rely on proxies. One such proxy is the Kullback-Leibler divergence, denoted by $D_{KL}(P||Q)$ and defined for discrete measures as
\begin{eqnarray*}
    D_{KL}(P||Q)&=&\sum_{x}P(x)\log\left(\frac{P(x)}{Q(x)}\right).
\end{eqnarray*}
The Kullback-Leibler divergence is related to total variation distance through the so-called Pinsker inequality. We refer the reader to \cite{Tsybakov} for a detailed discussion on this fundamental inequality.

\begin{Lemma}[Pinsker's Inequality]
    \label{pinsker}
    If $P$ and $Q$ are two probability measures defined on a measurable space $(\varOmega,\mathscr{F})$, then
    \begin{eqnarray*}
        \mathrm{TV}(P,Q) \leq\sqrt{\frac{D_{\mathrm{KL}}(P||Q)}{2}}.
    \end{eqnarray*}
\end{Lemma}

Pinsker's inequality implies that it is enough to show that $D_{\mathrm{KL}}(P||Q)$ is small to have a small total variation distance between $P$ and $Q$. Another consequence of this inequality is that convergence in total variation can be established through convergence in divergence. Specifically, we have that a sequence of probability measures $\{P_{n} : n\in\mathbb{N}\}$ converges in total variation distance to a probability measure $P$ if $\displaystyle \lim_{n\to\infty} D_{\mathrm{KL}}(P_{n}||P) = 0$. It is important to note that there exist significant relationships between the total variation distance and other important divergences, such as the Hellinger divergence or $\chi^{2}$-divergence \citep{polyanskiy2022information}. We focus on the Kullback-Leibler divergence due to its compatibility with the probabilistic structure of stick-breaking processes.

\section{Kullback-Leibler divergence for the Dirichlet and the geometric process}
\label{KL_DIR_GEOM}

In this section, and as a first approach to the study of the Kullback-Leibler divergence between stick-breaking priors, we centre our attention on the divergence between two random probability measures with common base measure, $P$ and $P'$, where $P$ is distributed as a Dirichlet process with iid length variables $\lengthVar\sim Be(1,\theta)$, and $P'$ is distributed as a geometric process with length variable $v\sim Be(a,b)$. Furthermore, for these two random measures, we assume that they share the locations of the jumps, $\locations$, so that the Kullback-Leibler divergence of $P$ with respect to $P'$ can be obtained as
\begin{eqnarray}
\label{KL_P_P'}
D_{\mathrm{KL}}(P||P')&=&\sum_{n=1}^{\infty}w_n\log\left(\frac{w_n}{w_n'}\right),
\end{eqnarray}
where $w_n$ are the weights of the Dirichlet process as in \eqref{SBP}, $w_n'$ the weights of the geometric process, that is, $w_n'=v(1-v)^{n-1}$, and where \eqref{KL_P_P'} holds with probability 1.

Moreover, since the Kullback-Leibler divergence is non-symmetric we can also obtain the divergence of $P'$ with respect to $P$ as
\begin{eqnarray*}
\label{KL_P'_P}
D_{\mathrm{KL}}(P'||P)&=&\sum_{n=1}^{\infty}w_n'\log\left(\frac{w_n'}{w_n}\right).
\end{eqnarray*}
Then as stated in Theorem~\ref{theo_KL}, a nice expression in terms of the length variables can be obtained for both divergences.

\begin{Theo}
\label{theo_KL}
Let $P$ be distributed as a Dirichlet process and $P'$ be distributed as a geometric process $(a,b,P_0)$. Then the Kullback-Leibler divergence of $P$ with respect to $P'$ is a.s. given by
\begin{eqnarray}
\label{Final_KL}
D_{\mathrm{KL}}(P||P')&=&\sum_{n=1}^{\infty}\left[\prod_{j=1}^{n-1}(1-v_j)\right]d(v_n||v),
\end{eqnarray}
where $d(v_n||v)$ is the binary divergence of $v_n$ with respect to $v$, {\it i.e.,},
\begin{eqnarray*}
d(v_n||v)&=&v_n\log\left(\frac{v_n}{v}\right)+(1-v_n)\log\left(\frac{1-v_n}{1-v}\right).
\end{eqnarray*}
Similarly, the Kullback-Leibler divergence of $P'$ with respect to $P$ is a.s. given by
\begin{eqnarray}
\label{Final_KL_2}
D_{\mathrm{KL}}(P'||P)&=&\sum_{n=1}^{\infty}(1-v)^{n-1}d(v||v_n).
\end{eqnarray}
\end{Theo}

For the proof of the first part of Theorem~\ref{theo_KL}, which we provide in Appendix~\ref{proofsSec3} of the supplementary material, we just require to note that, for a proper species sampling process, the sum of the tails of the weights satisfy the following relationship
\begin{eqnarray*}
\sum_{n=k+1}^{\infty}w_n&\stackrel{\text{a.s.}}{=}&\prod_{j=1}^k(1-v_j),
\end{eqnarray*}
which proof is straightforward. 

Since we are dealing with random measures the Kullback-Leibler divergence is a random variable as well; therefore, statistical properties such as the expectation and the variance can be obtained to have a better understanding of the divergence. To this end we have considered two different approaches depending on whether both processes have the same length variable starting point, that is, $v_1=v$ which we name  the coupled case and the uncoupled case otherwise. As we see in the following sections this is certainly an interesting way of tackling the problem. As a final note, it is important to remark that in the following sections we  provide the analysis of the random Kullback-Leibler divergence of the Dirichlet process with respect to the geometric process. This is due to mathematical tractability, ease of interpretation of the results derived and the connection made with Section~\ref{exch_sbp}, where we use the geometric process as the baseline for studying the divergence of the more general class of exchangeable stick-breaking processes. However, for completeness the reader can find in the Appendix~\ref{reversed} some interesting details about the expectation of the reversed Kullback-Leibler divergence.

\subsection{Uncoupled scenario}

In this section we are interested on the uncoupled scenario, that is, the case where $v_1\neq v$. To derive closed expressions for the expectation and the variance of the Kullback-Leibler divergence, it suffices to recall that for the Dirichlet process we have that $\lengthVar$ is a sequence of iid random variables with common $Be(1,\theta)$ distribution, whereas for the geometric process $v$ is a $Be(a,b)$ random variable. With this in mind and by taking into account the formulas of the expectation of logarithmically transformed beta random variables, which for ease we provide in Appendix~\ref{auxiliary_formulas} of the supplementary material, we are then able to present the expectation and the variance of the Kullback-Leibler divergence in Lemma~\ref{expectation_KL}.

\begin{Lemma}
\label{expectation_KL}
Let $P$ be distributed as a Dirichlet process and $P'$ be distributed as a geometric process, respectively defined by their stick-breaking length variables $\lengthVar\sim$ Be$(1,\theta)$ and $v\sim$ Be$(a,b)$. Then the expected value of the Kullback-Leibler divergence of $P$ with respect to $P'$ is given by
\begin{eqnarray*}
\Esp(D_{\mathrm{KL}}(P||P'))&=&(\theta+1)\Esp(d(v_1||v)),
\end{eqnarray*}
where 
\begin{eqnarray*}
\label{exp_div_v1_v}
\Esp(d(v_1||v))&=&\rho[\psi(2)-\psi(\theta+2)-\psi(a)+\psi(a+b)]+(1-\rho)[-\rho-\psi(b)+\psi(a+b)],
\end{eqnarray*}
and where $\rho=\frac{1}{\theta+1}$ is the tie probability of the Dirichlet process and $\psi(\cdot)$ is the digamma function. Finally, the corresponding variance is given by 
\begin{eqnarray*}
    \Var(D_{\mathrm{KL}}(P||P'))&=&(\theta+1)(\theta+2)\Esp((1-v_1)d(v_1||v)d(v_2||v))\nonumber\\
    &+&\frac{(\theta+2)}{2}\Esp(d^2(v_1||v))-[(\theta+1)\Esp(d(v_1||v))]^2.
\end{eqnarray*}
\end{Lemma}

\begin{figure}[ht!]
\centering
  \includegraphics[scale=.7]{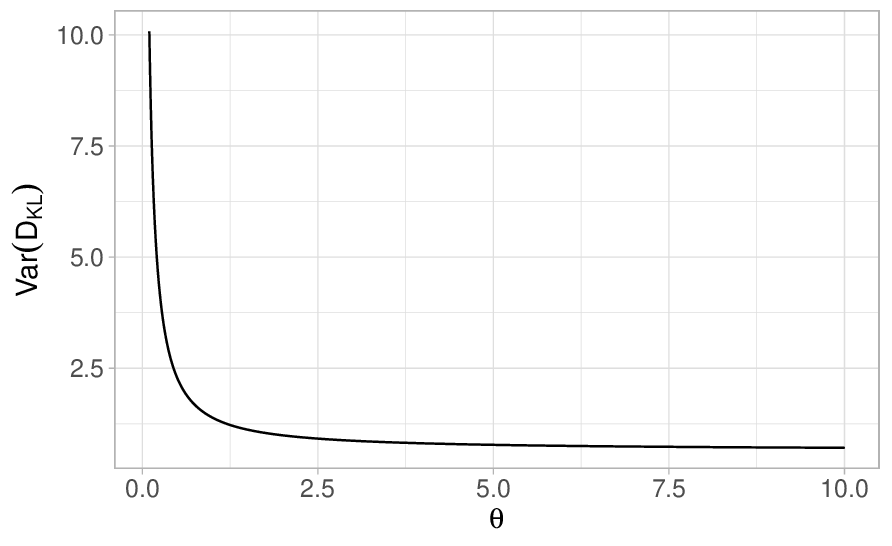}
\caption{Variance of the Kullback-Leibler divergence of the Dirichlet process $(\theta)$ with respect to the geometric process $(1,\theta)$ as a function of $\theta$.}
\label{var_uncoupled}
\end{figure}

The proof of Lemma~\ref{expectation_KL} can be seen in Appendix~\ref{proofsSec3} of the supplementary material, where we also provide closed expressions for $\Esp(d^2(v_1||v))$ and $\Esp((1-v_1)d(v_1||v)d(v_2||v))$. Certainly, Lemma~\ref{expectation_KL} provides us with the necessary tools to study the behaviour of the divergence as a function of $\theta,a$ and $b$. Finally, it is definitely interesting to notice that whenever $a=1$ and $b=\theta$ (which is a case we are  interested in this paper), the expectation of the Kullback-Leibler divergence does not longer depend on the value of $\theta$, in fact after some calculations, it is straightforward to see that $\Esp(D_{\mathrm{KL}}(P||P'))=1$. So that obtaining and analysing the variance is crucial to understand the variability as a function of $\theta$. For this, in Figure~\ref{var_uncoupled} we plot the variance as a function of $\theta$ for which we are able to see a monotonic decreasing behaviour that approximately converges to 0.6449. This is certainly a compelling result; however, since we are assuming both processes share the set of locations, $\locations$, we can then analyse the coupled case, that is, $v_1=v$, so that both processes start at the same length variable starting point and a more direct comparison can be derived.

\subsection{Coupled scenario}
As briefly explained in the previous section, an interesting case happens when $v_1=v$ which we have denoted as the coupled scenario. As we shall see in Section~\ref{exch_sbp}, this is something that allows us to fully recover the geometric process by tuning the tie probability in a similar fashion as in \cite{Leyva_2021}. Now, first we can analyse the statistical properties of the Kullback-Leibler divergence of the Dirichlet and the geometric process for the case $v_1=v$. In this situation we are also able to provide closed expressions for the expectation and the variance as detailed in Lemma~\ref{expectation_coupled_KL}.

\begin{Lemma}
\label{expectation_coupled_KL}
Let $P$ be distributed as a Dirichlet process and $P'$ be distributed as a geometric process, respectively defined by their stick-breaking length variables $\lengthVar\sim$ Be$(1,\theta)$ and $v=v_1$. Then the expected value of the Kullback-Leibler divergence of $P$ with respect to $P'$ is given by,
    \begin{eqnarray*}
    \Esp(D_{\mathrm{KL}}(P||P'))&=&\frac{\theta}{\theta+1}.
    \end{eqnarray*}
    Moreover, the variance is given by
    \begin{eqnarray*}
    \Var(D_{\mathrm{KL}}(P||P'))&=&(\theta+2)(\theta+1)\Esp((1-v_1)^2(1-v_2)d(v_2||v_1)d(v_3||v_1))\nonumber\\&+&\frac{(\theta+2)}{2}\Esp(((1-v_1)d(v_2||v_1))^2)-(1-\rho)^2.
\end{eqnarray*}
\end{Lemma}

The proof of Lemma~\ref{expectation_coupled_KL} can be seen in Appendix~\ref{proofsSec3} of the supplementary material, where we also provide closed expressions for $\Esp(((1-v_1)d(v_2||v_1))^2)$ and $\Esp((1-v_1)^2(1-v_2)d(v_2||v_1)d(v_3||v_1))$ which again depend on the expectation of logarithmically transformed beta random variables. For the coupled scenario, we are able to see that the expectation has a nice expression in terms of the tie probability of the Dirichlet process since $\theta(\theta+1)^{-1}=1-\rho$. 

\begin{figure}[ht!]
\centering
  \includegraphics[scale=.7]{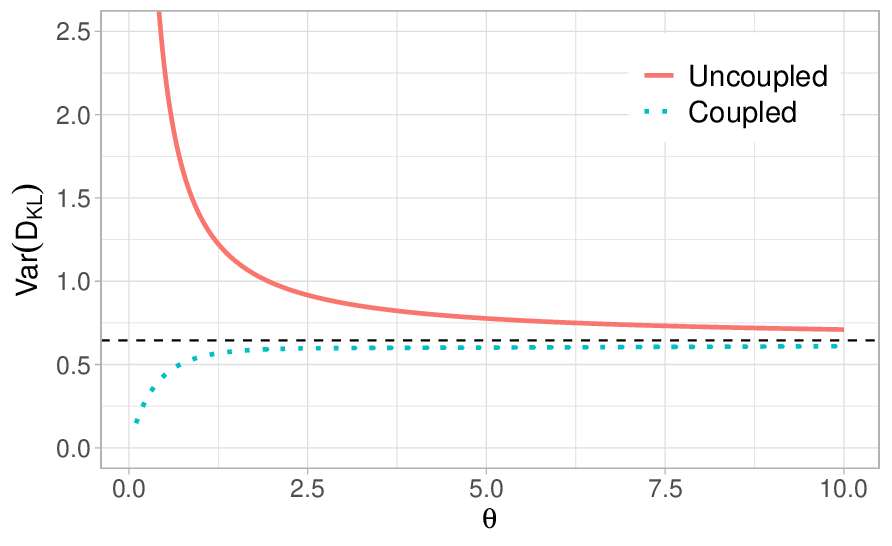}
\caption{Comparison of the variances of the Kullback-Leibler divergence of the uncoupled (red line) and coupled Dirichlet and geometric process (blue dots), and where the black dashed line represents the limit as $\theta\to\infty$ of both variances.}
\label{var_comps}
\end{figure}

Now, there are some interesting remarks that can be directly seen from the expectation of the Kullback-Leibler divergence. Firstly, $\Esp(D_{\mathrm{KL}}(P||P'))\to1$ as $\theta\to\infty$ which is the expectation of the divergence for the uncoupled case with $a=1$ and $b=\theta$. Secondly, as $\theta\to0$, we have that $\Esp(D_{\mathrm{KL}}(P||P'))\to0$, something that was clearly expected, since it corresponds to a degenerate case where all the stick-breaking length variables are a.s. one. Finally, as to the variance of the Kullback-Leibler divergence as a function of $\theta$, we can observe a monotonic increasing behaviour contrary to the uncoupled case as illustrated in Figure~\ref{var_comps}. However, it is certainly quite compelling to remark that in this case we also get a convergence to (approximately) 0.6449. 

\section{KL divergence for stick-breaking processes with \linebreak exchangeable length variables}
\label{exch_sbp}
Now that we have provided a closed expression for the Kullback-Leibler divergence of the Dirichlet process with respect the geometric process and viceversa, we can turn our attention to the more general case of stick breaking processes with exchangeable length variables. This should allow us to connect the theory developed so far with the results presented in \cite{Leyva_2021}. In particular we are interested on characterising the expected value of the Kullback-Leibler divergence in terms of the tie probability of the length variables, {\it i.e.,} $\rho_{\nu}$, from which we know, that under certain conditions we can then recover the Dirichlet and the geometric processes as limit cases. For the purposes of this paper we restrict our attention to Dirichlet-driven length variables, that is, 
\begin{eqnarray*}
 v_i|\nu&\simiid&\nu\\
    \nu&\sim&\text{DP}(\beta,\nu_0)\nonumber,
\end{eqnarray*}
where $\beta$ is the concentration parameter and $\nu_0$ is the base measure.

\subsection{Dirichlet-driven length variables}
As a first approach to the theory of exchangeable length variables, we consider the Dirichlet process as the nonparametric prior due to its inherent flexibility and mathematical tractability. To this end, we can recall that the Dirichlet process is completely specified by its base measure and its concentration parameter. Since we are dealing with a sequence of random variables defined in $\UnitInt$, we let the base measure $\nu_0=Be(1,\theta)$, which further allows us to guarantee that we have a proper species sampling model \citep{Leyva_2021} and denote by $\beta$ the concentration parameter, so that the tie probability is given by $\rho_{\nu}=(\beta+1)^{-1}$. Furthermore, the EPPF $p(\cdot)$ is given \citep[\it{see e.g.} ][]{EWENS1972,Antoniak1974} by 
\begin{eqnarray}
\label{eppf}
    p(n_1,\ldots,n_k)&=&\frac{\beta^k}{(\beta)^{(n)}}\prod_{i=1}^k(n_i-1)!,
\end{eqnarray}
where $(\beta)^{(n)}$ is the Pochhammer symbol, {\it i.e.,}, $(\beta)^{(n)}=\beta(\beta+1)\cdots(\beta+n-1)$. 

We are then able to recover the expression of the Kullback-Leibler divergence of a stick-breaking process with exchangeable sequence of length variables $\lengthVar$ driven by a Dirichlet process, henceforth denoted by $P_{\beta}$, with respect to the geometric process as stated in Theorem~\ref{KL_SBP_GEOM} and whose proof can also be seen in Appendix~\ref{proofsSec4} of the supplementary material.
\begin{Theo}
\label{KL_SBP_GEOM}
Let $P_{\beta}$ a stick-breaking process with exchangeable length variables driven by a Dirichlet process $(\beta,\nu_0)$, with $\nu_0=Be(1,\theta)$ and $P'$ a random measure distributed as a geometric process $(1,\theta,P_0)$. Then the Kullback-Leibler divergence of $P_{\beta}$ with respect to $P'$ is a.s. given by
\begin{eqnarray*}
\label{kl_p_beta_p}
D_{\mathrm{KL}}(P_{\beta}||P')&=&\sum_{n=1}^{\infty}\left[\prod_{k=1}^{n-1}(1-v_k)\right]d(v_n||v),
\end{eqnarray*}
where $d(v_n||v)$ is the binary divergence of $v_n$ with respect to $v$.
\end{Theo}

Then using the EPPF of the Dirichlet process \eqref{eppf} and the formula described in \cite{Leyva_2021} for the expectation of a measurable function $f$ of exchangeable random variables, we are able to first provide an expression for the expectation of the Kullback-Leibler divergence as stated in Lemma~\ref{exp_exch_sbp}, whose proof we provide in Section B of the supplementary material.

\begin{Lemma}
    \label{exp_exch_sbp}
    Let $P_{\beta}$ and $P'$ as in Theorem~\ref{KL_SBP_GEOM}. Then the expectation the Kullback-Leibler divergence of $P_{\beta}$ with respect to $P'$ is a.s. given by
\begin{eqnarray*}
\label{exp_kl_p_beta_p}
    D_{\theta}(\beta)&=&\Esp[D_{\mathrm{KL}}(P_{\beta}||P')]=\sum_{n=1}^{\infty}\sum_{\pi\in\mathscr{P}([n])}F_{\theta}(\pi)p_{\beta}(\pi),
\end{eqnarray*}
where
\begin{eqnarray*}
F_{\theta}(\pi)&=&\Esp\left(\left(\prod_{j=1}^{k-1}(1-v_j)^{n_j}\right)(1-v_{k})^{n_k-1}d(v_k||v)\right),
\end{eqnarray*}
and
\begin{eqnarray*}
p_{\beta}(\pi)&=&\frac{\beta^{k}}{(\beta)^{(n)}}\prod_{j=1}^{k}(n_j-1)!
\end{eqnarray*}
\end{Lemma}

It is then quite evident that obtaining a closed expression is mathematically challenging since there is a need to sum over all the partitions. However, this expression provides us with a nice way to connect the theory developed so far, with the idea of recovering the dependent and independent length-variables through the means of tuning $\beta$. In this direction, we can first notice that for $\beta=0$, we will have by definition of the Dirichlet-driven exchangeable length variables that $v_1=v_2=\cdots$, so that $\forall n\geq1$ $p_{0}(\pi)=\delta_{(1_n)}(\pi)$,
where $\delta$ is the Dirac measure and $1_n\in\mathscr{P}([n])$ is the partition where all the variables are the same. Hence,
\begin{eqnarray*}
        D_{\theta}(0)&=&\sum_{n=1}^{\infty}\Esp[(1-v_1)^{n-1}d(v_1||v)],
\end{eqnarray*}
which further allows us to recover some of the expressions derived in Section~\ref{KL_DIR_GEOM}. In particular we know that for $v_1\neq v$, which corresponds to the uncoupled case, we obtain that $D_{\theta}(0)=\infty$, and for the coupled case, that is $v_1=v$, we have that $D_{\theta}(0)=0$. 

\begin{figure}[ht!]
\centering
  \includegraphics[scale=.7]{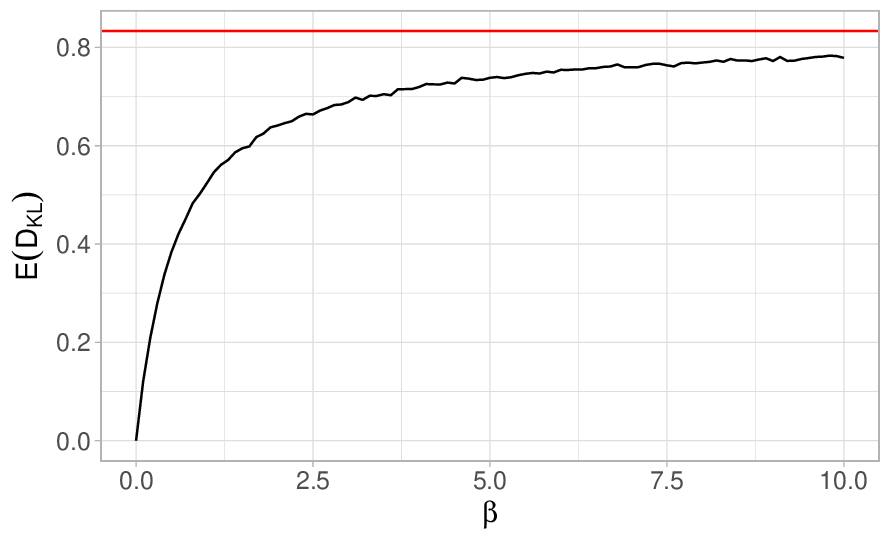}
\caption{Expectation of the Kullback-Leibler divergence of $P_{\beta}$ with respect to $P$ for different values of $\beta$ where the red line is equal to $\theta(\theta+1)^{-1}$ for $\theta=5$.}
\label{exp_coupled_KL}
\end{figure}

Intuitively, for the case where we take $\beta\to\infty$, we should be able to recover the expectation of the Kullback-Leibler of the Dirichlet process with respect to the geometric one, since this would imply that all the variables are iid drawn from the common $\text{Be}(1,\theta)$ distribution. This is something that can be appreciated in Figure~\ref{exp_coupled_KL}, where for $n=300$ length variables we obtain the $D_{\mathrm{KL}}(P_{\beta}||P')$ for several values of $\beta\in[0,10]$ and $\theta=5$, so that the red line represents the expectation of the Kullback-Leibler divergence of the Dirichlet process with respect to the geometric one and given by $\theta(\theta+1)^{-1}=0.8333$. Of course, since closed-form expressions are not at hand, we are recurring to the classical Monte Carlo estimate by sampling and averaging 100,000 simulations for each scenario of $\beta$ considered.

At this point, it becomes quite important to remark that analysing the properties of $D_{\theta}(\beta)$ is certainly not straightforward and it is mathematically challenging, which is something quite common when dealing with these Bayesian nonparametric models, where there is a need to consider the sum over all the possible partitions. For example, at first sight and based on our intuition and knowledge on the behaviour of the exchangeable stick-breaking processes, we conjecture that the expectation is monotonically increasing with respect to $\beta$ for fixed $\theta$. Of course, proving this statement is certainly not trivial. Other properties such as continuity and boundedness for specific values of $\beta$ and $\theta$ are addressed in Theorem~\ref{bound_cont_comp} which allows us to see that, for $\beta>\theta+1$, we have that $D_{\theta}(\beta)$ is a continuous function for which we are able to further provide an upper bound. 

\begin{Theo}
\label{bound_cont_comp}
Let $P_{\beta}$ a stick-breaking process with base measure $P_0$ and with an exchangeable length variables $\lengthVar$ driven by a Dirichlet process $(\beta,\nu_0)$ where $\nu_0=Be(1,\theta)$, $P$ a Dirichlet process $(\theta,P_0)$, and $P'$ a geometric process $(1,\theta,P_0)$. Then, for $\beta>\theta+1$, we have that $\beta \mapsto D_{\theta}(\beta)$ is continuous,
\begin{eqnarray*}
    \label{upper_bound}
    D_{\theta}(\beta)\leq\frac{\theta}{\theta+1}\frac{\beta^2}{(\beta-1)(\beta-(\theta+1)},
\end{eqnarray*}
and, in addition, 
\begin{eqnarray*}
\lim_{\beta\rightarrow\infty}\Esp(D_{\mathrm{KL}}(P_{\beta}||P'))&=&\Esp(D_{\mathrm{KL}}(P||P')).
\end{eqnarray*}
\end{Theo}

The proof of Theorem~\ref{bound_cont_comp} is provided in Appendix~\ref{proofsSec4} of the supplementary material. It relies on Lemma~\ref{Lemma:SumRisingFactorials} and Lemma~\ref{Lemma:SumQuotientRisingFactorials} regarding the rising factorial. The proofs of the latter lemmas rely on the theory of reciprocals of inverse factorial series \cite[\textit{see e.g.}][]{norlund, harris_jr1964,tweddle2012james}, and that can also be found in Appendix~\ref{proofsSec4} of the supplementary material.

\begin{Lemma}
\label{Lemma:SumRisingFactorials}
If $\beta > 1$, then 
\begin{eqnarray*}
\sum_{n=1}^{\infty} \frac{(n-1)!}{\beta^{(n)}} &=& \frac{1}{\beta-1}.
\end{eqnarray*}
\end{Lemma}

\begin{Lemma}
\label{Lemma:SumQuotientRisingFactorials}
Let $\lambda\in(0,1)$ and $\beta>0$. If $(1-\lambda)\beta>1$, then 
\begin{eqnarray*}
    \sum_{n=1}^{\infty} \frac{(\lambda\beta)^{(n)}}{\beta^{(n)}} &=& \frac{\lambda\beta}{(1-\lambda)\beta-1}.
\end{eqnarray*}
\end{Lemma}

\begin{figure}[ht!]
\centering
  \includegraphics[scale=.7]{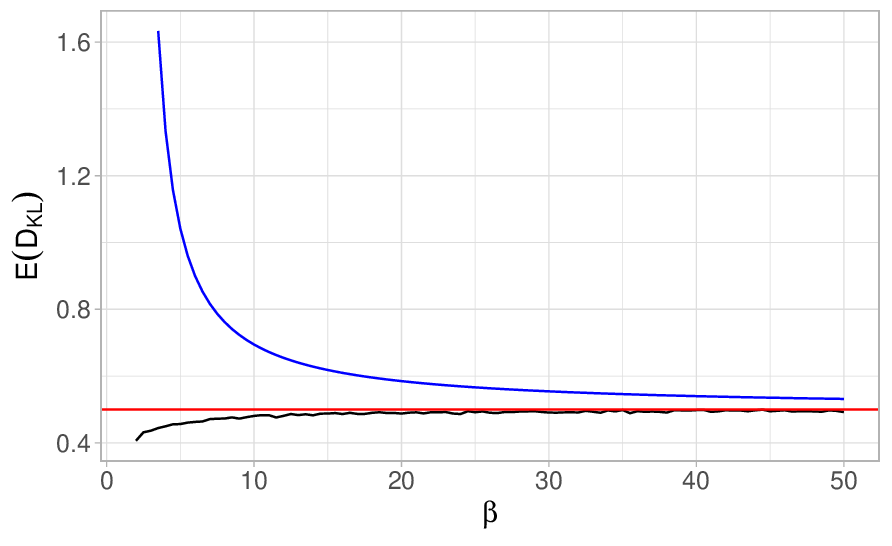}
\caption{Expectation of the Kullback-Leibler divergence of $P_{\beta}$ with respect to $P$ (black line) for $\beta\in(\theta+1,50)$ and $\theta=1$, where the red line represents the limiting behaviour of the expectation as $\beta\to\infty$, and where the blue line represents the upper bound on the expectation.}
\label{bound_figure}
\end{figure}

Note that we can fully recover the results of Section~\ref{KL_DIR_GEOM}. Specifically, we have that as $\beta\to\infty$ we recover the expectation of the Kullback-Leibler divergence of the Dirichlet process with respect to the geometric one, which is given by $\theta(\theta+1)^{-1}$ in the coupled case, that is, $v_1=v$. To illustrate the content of Theorem~\ref{bound_cont_comp}, we direct the attention of the reader to Figure~\ref{bound_figure}, where it can be appreciated the limiting behaviour of the expectation of the random Kullback-Leibler divergence, as well as the upper bound which gets tighter as $\beta$ increases. Although Theorem~\ref{bound_cont_comp} only covers the case $\beta > \theta+1$, exhaustive empirical evidence and related theoretical results suggest that the same is true for all $\beta>0$. As a result, we put forward the following conjecture.

\begin{Conj}
\label{Conjecture}
Let $P_{\beta}$ a stick-breaking process with base measure $P_0$ and with an exchangeable length variables $\lengthVar$ driven by a Dirichlet process $(\beta,\nu_0)$ where $\nu_0=Be(1,\theta)$, $P$ a Dirichlet process $(\theta,P_0)$, and $P'$ a geometric process $(1,\theta,P_0)$. Then, for $\beta>0$, we have that $\beta \mapsto D_{\theta}(\beta)$ is continuous and $\displaystyle D_{\theta}(\beta)\leq\frac{\theta}{\theta+1}$.
\end{Conj}

Recall Pinsker's inequality as presented in Section~\ref{Subsection:TotalVariationDistance}. Jensen's inequality and the previous conjecture imply that
\begin{eqnarray*}
    \Esp[\mathrm{TV}(P_{\beta},P')]\leq\sqrt{\frac{\Esp[D_{\mathrm{KL}}(P_{\beta}||P')]}{2}} \leq \sqrt{\frac{\theta}{2(\theta+1)}}.
\end{eqnarray*}
Hence, when $\theta$ is small, the complexity cost of introducing $\beta$ does not pay-off in terms of the flexibility of the model, as all prior distributions $P_{\beta}$ are close in total variation distance.

Also, the continuity posed in Conjecture~\ref{Conjecture} implies that
\begin{eqnarray*}
    \Esp[\mathrm{TV}(P_{\beta},P')]&\leq&\sqrt{\frac{\Esp[D_{\mathrm{KL}}(P_{\beta}||P')]}{2}}\stackrel{\beta\to 0}{\longrightarrow}0.
\end{eqnarray*}
This convergence in total variation distance recovers, in particular, the weak convergence of $P_{\beta}$ to the geometric process demonstrated in \cite{Leyva_2021}.

\section{Discussion}
\label{conclusions}
In this paper we have centred our attention on a methodology to analyse the discrepancies between extreme classes of stick-breaking process through the means of the Kullback-Leibler divergence. This information-based prior analysis has allowed us to derive interesting analytical results about the divergence of the Dirichlet process and the geometric process, specifically we were able to provide a mean-variance characterisation of this random divergence. Centring our attention on the divergence of the Dirichlet process with respect to the geometric process, an appealing framework to study the divergence of the more general exchangeable stick-breaking processes emerged. This is particularly useful, since our results allows us to provide conditions under which all of the stick-breaking processes considered exhibit close total variation. So that, for $\beta$ large the complexity cost of introducing the exchangeable stick-breaking process does not pay-off in terms of performance for small $\theta$. This is is of course very important when eliciting a nonparametric prior.

\section*{Acknowledgements}
We dedicate this paper to the memory of our co-author, Dr. Mario Diaz, who passed away during the final stages of this project. His insight, dedication, and enthusiasm were fundamental to the development of this work, and his absence is profoundly felt. We are grateful to the editor and an anonymous reviewer for their valuable comments that help us to improve the manuscript.

\section*{Disclosure statement}
The authors report there are no competing interests to declare.

\section*{Funding}
The authors gratefully acknowledge the support of PAPIIT Grant IT100524. 

\section*{ORCID}
José A. Perusquía \url{https://orcid.org/0000-0001-7722-7711}\\
Mario Diaz \url{https://orcid.org/0000-0002-9321-9815}\\
Ramsés H. Mena \url{https://orcid.org/0000-0001-9608-8059}

\bibliographystyle{hapalike}
\bibliography{bibliography.bib}
\newpage
\section*{Supplementary material}
\begin{appendix}
\section{Expectation of logarithmically transformed beta \linebreak random variables}
\label{auxiliary_formulas}
In this section, we present some already known and useful formulas, regarding the expectation of logarithmically transformed beta random variables, which are used to obtain the expectation of the Kullback-Leibler divergence of the Dirichlet process and the geometric process. To this end, if we let $X\sim Be(a,b)$ then
\begin{eqnarray*}
\Esp(\log(X))&=&\psi(a)-\psi(a+b),\\
\Esp(\log(1-X))&=&\psi(b)-\psi(a+b),\\
\Esp(\log^2(X))&=&(\psi(a)-\psi(a+b))^2+(\psi_1(a)-\psi_1(a+b)),\\
\Esp(\log^2(1-X))&=&(\psi(b)-\psi(a+b))^2+(\psi_1(b)-\psi_1(a+b)),\\
\Esp(\log(X)\log(1-X))&=&(\psi(a)-\psi(a+b))(\psi(b)-\psi(a+b))-\psi_1(a+b),
\end{eqnarray*}
where $\psi(\cdot)$ is the digamma function defined as,
\begin{eqnarray*}
    \psi(x)&=&\frac{\partial \log[\varGamma(x)]}{\partial x} = \frac{\varGamma(x)'}{\varGamma(x)},
\end{eqnarray*}
and $\psi_1(\cdot)$, is the trigamma function defined as the first derivative of $\psi(\cdot)$. Hence, it is straightforward to see that, for suitable values $c$ and $d$,
\begin{eqnarray*}
\Esp(X^c(1-X)^d\log(X))&=& \frac{B(a+c,b+d)}{B(a,b)} \Esp(\log(X')),\\
\Esp(X^c(1-X)^d\log(1-X))&=&\frac{B(a+c,b+d)}{B(a,b)}\Esp(\log(1-X')),\\
\Esp(X^c(1-X)^d\log^2(X))&=&\frac{B(a+c,b+d)}{B(a,b)}\Esp(\log^2(X')),\\
\Esp(X^c(1-X)^d\log^2(1-X))&=&\frac{B(a+c,b+d)}{B(a,b)}\Esp(\log^2(1-X')),\\
\Esp(X^c(1-X)^d\log(X)\log(1-X))&=&\frac{B(a+c,b+d)}{B(a,b)}\Esp(\log(X')\log(1-X')),
\end{eqnarray*}
where $X'\sim Be(a+c,b+d)$.

\clearpage
\section{Reversed divergence}
\label{reversed}
Now that we have analysed the Kullback-Leibler divergence of the Dirichlet process with respect to the geometric for both uncoupled and coupled cases, we direct our attention to the reversed divergence. To this end we recall that in Theorem~\ref{theo_KL} we have already derived an expression for the divergence as stated in \eqref{Final_KL_2}; however, as we see from Lemma~\ref{reversed_exp}, in this case the expectation becomes quite more challenging and a closed expression is not as straightforward as before, although it is important to remark that it can be obtained numerically using specialised software.

\begin{Lemma}
\label{reversed_exp}
Let $P$ be distributed as a Dirichlet process and $P'$ be distributed as a geometric process, respectively defined by their stick-breaking length variables $\lengthVar\sim$ Be$(1,\theta)$ and $v\sim$ Be$(a,b)$. Then the expected value of the Kullback-Leibler divergence of $P'$ with respect to $P$ is given by,

\begin{eqnarray*}
\Esp(D_{\mathrm{KL}}(P'||P))\nonumber&=&\sum_{n=1}^{\infty}\biggl(\frac{a\varGamma(b+n-1)\varGamma(a+b)}{\varGamma(b)\varGamma(a+b+n)}\left[\psi(a+1)-\psi(a+b+n)+\gamma+\psi(\theta+1)\right]\nonumber\\    &&+\frac{\varGamma(b+n)\varGamma(a+b)}{\varGamma(b)\varGamma(a+b+n)}\left[\psi(b+n)-\psi(a+b+n)+\frac{1}{\theta}\right]\biggr),
\end{eqnarray*}
where $\gamma=-\psi(1)$ is the Euler-Mascheroni constant, and $\Esp(D_{\mathrm{KL}}(P'||P))<\infty$ only if $a>1$.
\end{Lemma}
\begin{proof}[Proof of Lemma~\ref{reversed_exp}]
Consider the Kullback-Leibler divergence of the geometric process with respect the Dirichlet process, which is given by,
\begin{eqnarray*}
    D_{\mathrm{KL}}(P'||P)&=&\sum_{n=1}^{\infty}(1-v)^{n-1}d(v||v_n).
\end{eqnarray*}
By Beppo-Levi's theorem, we have that,
\begin{eqnarray*}
    \Esp(D_{\mathrm{KL}}(P'||P))&=&\sum_{n=1}^{\infty}\Esp\left((1-v)^{n-1}d(v||v_n)\right)\\
    &=&\sum_{n=1}^{\infty}\Esp\left((1-v)^{n-1}\left(v\log\left(\frac{v}{v_n}\right)+(1-v)\log\left(\frac{1-v}{1-v_n}\right)\right)\right).
\end{eqnarray*}
This expectation can be split into the following four terms,
\begin{align*}
f_1&=\Esp\left[v(1-v)^{n-1}\log(v)\right],\\
f_2&=\Esp\left[v(1-v)^{n-1}\right]\Esp\left[\log(v_n)\right],\\
f_3&=\Esp\left[(1-v)^{n}\log(1-v)\right],\\
f_4&=\Esp\left[(1-v)^{n}\right]\Esp\left[\log(1-v_n)\right].
\end{align*}
Then, using the same reasoning as for Lemma~\ref{expectation_KL} and the formulas provided in Appendix~\ref{auxiliary_formulas}, it is straightforward to see that,
\begin{align*}
f_1&=\frac{B(a+1,b+n-1)}{B(a,b)}[\psi(a+1)-\psi(a+b+n)],\\
f_2&=\frac{B(a+1,b+n-1)}{B(a,b)}[\psi(1)-\psi(\theta+1)],\\
f_3&=\frac{B(a,b+n)}{B(a,b)}[\psi(b+n)-\psi(a+b+n)],\\
f_4&=\frac{B(a,b+n)}{B(a,b)}[\psi(\theta)-\psi(\theta+1)],
\end{align*}
where 
\begin{align*}
    \frac{B(a+1,b+n-1)}{B(a,b)}&=\frac{a\varGamma(b+n-1)\varGamma(a+b)}{\varGamma(a+b+n)\varGamma(b)},\\
    \frac{B(a,b+n)}{B(a,b)}&=\frac{\varGamma(b+n)\varGamma(a+b)}{\varGamma(a+b+n)\varGamma(b)}.
\end{align*}
Hence, by joining these four expressions we get the desired result about the expectation of the divergence. Now, if we let $a=1+p$, with $p\in(0,1)$, and using Gautschi's inequality we obtain that,
\begin{eqnarray*}
    \frac{a\varGamma(b+n-1)\varGamma(a+b)}{\varGamma(a+b+n)\varGamma(b)}&=&\frac{(1+p)\varGamma(b+n-1)\varGamma(b+p+1)}{\varGamma(b+n+p+1)\varGamma(b)}\\
    &=&\frac{(1+p)\varGamma(b+n+1)\varGamma(b+p+1)}{(b+n-1)(b+n)(b+n+p)\varGamma(b+n+p)\varGamma(b)}\\
    &<&\frac{(1+p)(b+n+1)^{1-p}\varGamma(b+p+1)}{(b+n-1)(b+n)\varGamma(b)}\\
    &\sim& \left(\frac{1}{n^{1+p}}\right).
\end{eqnarray*}
Similarly, for the second coefficient we have that
\begin{eqnarray*}
    \frac{\varGamma(b+n)\varGamma(a+b)}{\varGamma(a+b+n)\varGamma(b)}&=&\frac{\varGamma(b+n)\varGamma(b+1+p)}{\varGamma(b+n+p+1)\varGamma(b)}\\
    &=&\frac{\varGamma(b+n+1)\varGamma(b+1+p)}{(b+n)(b+n+p)\varGamma(b+n+p)\varGamma(b)}\\
    &<&\frac{(b+n+1)^{1-p}\varGamma(b+1+p)}{(b+n)(b+n+p)\varGamma(b)}\\
    &\sim&\left(\frac{1}{n^{1+p}}\right).
\end{eqnarray*}
Note that such series are convergent for $p>0$. Now, for the digamma function we know that,
\begin{eqnarray*}
    \psi(a+b+n)&\sim& \log(n)-\frac{1}{2n},
\end{eqnarray*}
and that
\begin{eqnarray*}
    \psi(b+n)-\psi(a+b+n)&=& \psi(b+n)-\psi(1+p+b+n)\\
    &\leq&\psi(b+n)-\psi(1+b+n)\\
    &=&-\frac{1}{b+n}\\
    &\sim&-\frac{1}{n}.
\end{eqnarray*}
Hence, the four terms of the expectation are convergent only if $p>0$. Finally, as $p\to0$, we have that the coefficients are of order $n^{-1}$ and the divergence of the series is immediate.
\end{proof}

\clearpage
\section{Proofs of main results}
\subsection{Proofs of section 3}
\label{proofsSec3}
\begin{proof}[Proof of Theorem~\ref{theo_KL}]
Let $P$ be distributed as a Dirichlet process, and $P'$ be distributed as a geometric process as stated in the theorem. Observe that, at a formal level,
\begin{eqnarray}
    \label{DKL}
    D_{\mathrm{KL}}(P||P')&=&\sum_{n=1}^{\infty}w_n\log\left(\frac{w_n}{w_n'}\right)\nonumber\\
    &=&v_1\log\left(\frac{v_1}{v}\right)+\sum_{n=2}^{\infty}w_n\log\left(\frac{v_n}{v}\prod_{k=1}^{n-1}\left[\frac{1-v_k}{1-v}\right]\right)\nonumber\\
    \label{eq:DKDecomposition} &=&v_1\log\left(\frac{v_1}{v}\right)+\sum_{n=2}^{\infty}w_n\log\left(\frac{v_n}{v}\right)+\sum_{n=2}^{\infty}w_n\sum_{k=1}^{n-1}\log\left(\frac{1-v_k}{1-v}\right).\nonumber\\
\end{eqnarray}

Next we show that the previous formula indeed holds almost surely, by showing that the corresponding series converge absolutely almost surely. First we see that, by the strong law of large numbers, there exists a set $A_1$, such that $\mathbbm{P}(A_1) = 1$ and, for all $\omega \in A_1$,
\begin{equation}
\label{eq:AlmostSureConvergence1}
    \frac{1}{n} \sum_{k=1}^{n-1} \log(1-v_{k}) \to \mathbbm{E}[\log(1-V)].
\end{equation}
Then, let $\delta>0$ and observe that
\begin{equation*}
    \sum_{n=1}^{\infty} \mathbbm{P}\bigg(\bigg\lvert\frac{\log(v_{n})}{n}\bigg\rvert > \delta\bigg) = \sum_{n=1}^{\infty} 1-(1-\exp(-n\delta))^\theta<\infty,
\end{equation*}
where the last inequality follows directly from the d'Alambert's criterion since
\begin{equation*}
    \lim_{n\to\infty}\Big|\frac{1-(1-\exp(-((n+1)\delta))^\theta}{1-(1-\exp(-n\delta))^\theta}\Big|=\exp(-\delta)<1.
\end{equation*}
Hence, the Borel-Cantelli lemma implies that there exists a set $A_2$, such that $\Pro(A_2)=1$ and, for all $\omega\in A_2$,
\begin{equation}
\label{eq:AlmostSureConvergence2}
    \frac{\log(v_{n})}{n} \to 0.
\end{equation}

Now, if we define $A = A_1 \cap A_2$, we can observe that $\Pro(A)=1$, and in that set it is satisfied
\begin{equation*}
    \frac{\log(w_{n})}{n} = \frac{\log(v_{n}) + \sum_{k<n} \log(1-v_{k})}{n} \to \mathbbm{E}[\log(1-V)].
\end{equation*}
In particular, for all $\omega \in A$ and every $\epsilon>0$, there exists $N = N(\omega,\epsilon)$ such that, for all $n>N$,
\begin{equation*}
    \frac{\log(w_{n})}{n} \leq \mathbbm{E}[\log(1-V)] + \epsilon,
\end{equation*}
or, equivalently,
\begin{equation*}
    w_{n} \leq \big(e^{\mathbbm{E}[\log(1-V)]+\epsilon}\big)^{n}.
\end{equation*}
Taking $\epsilon < \lvert \mathbbm{E}[\log(1-V)] \rvert$ and defining $\lambda = e^{\mathbbm{E}[\log(1-V)]+\epsilon} < 1$, we conclude that, almost surely, for $n$ sufficiently large,
\begin{equation*}
    w_{n} \leq \lambda^{n}.
\end{equation*}

Returning to the absolute convergence of the series in \eqref{eq:DKDecomposition}, note that for the second series we have that
\begin{eqnarray*}
\sum_{n=2}^{\infty}\Big|w_n\log\left(\frac{v_n}{v}\right)\Big|&\leq&\sum_{n=2}^{\infty}w_n|\log(v_n)|+\sum_{n=2}^{\infty}w_n|\log(v)|\\
&=&-\sum_{n=2}^{\infty}w_n\log(v_n)-\sum_{n=2}^{\infty}w_n\log(v),
\end{eqnarray*}
where $-\sum_{n=2}^{\infty}w_n\log(v)\stackrel{\text{a.s.}}{\rightarrow}-\log(v)(1-w_1)$. Also, observe that the series
\begin{eqnarray*}
-\sum_{n=2}^{\infty}w_n\log(v_n)&=&\sum_{n=2}^{\infty}(n w_n) \frac{-\log(v_n)}{n},
\end{eqnarray*}
converges almost surely as, in the set $A$, the weights $w_{n}$ decrease geometrically for sufficiently large $n$, {\it i.e.,}, $w_n<\lambda^n$ with $\lambda<1$. Hence, we have that
\begin{eqnarray*}
\sum_{n=2}^{\infty}|n w_n|&<&\infty.
\end{eqnarray*}

Then, since $\frac{-\log(v_n)}{n}$ is a non-negative sequence, for which we know it converges to zero by \eqref{eq:AlmostSureConvergence2}, we obtain the a.s. convergence of $\sum_{n=2}^{\infty}w_n\log(v_n)$. Therefore, the series $\sum_{n\geq2}w_n\log\left(\frac{v_n}{v}\right)$ converges absolutely almost surely. Now, in a similar fashion we have that,
\begin{eqnarray*}
\sum_{n=2}^{\infty}\Big|w_n\sum_{k=1}^{n-1}\log\left(\frac{1-v_k}{1-v}\right)\Big|&\leq&\sum_{n=2}^{\infty}w_n\sum_{k=1}^{n-1}\Big|\log\left(\frac{1-v_k}{1-v}\right)\Big|\\
&\leq&\sum_{n=2}^{\infty}w_n\sum_{k=1}^{n-1}|\log(1-v_k)|+|\log(1-v)|\sum_{n=2}^{\infty}(n-1)w_n\\
&=&-\sum_{n=2}^{\infty}(nw_n)\sum_{k=1}^{n-1}\frac{\log(1-v_k)}{n}-\log(1-v)\sum_{n=2}^{\infty}(n-1)w_n.
\end{eqnarray*}
Note that the second series converge a.s. since $\sum_{n=2}^{\infty}(n-1)w_n$ converges a.s., because the weights $w_{n}$ decrease geometrically. The absolute convergence of the first series follows from the previous fact and \eqref{eq:AlmostSureConvergence1}. Thus, we have absolute convergence as desired. Now, since we have absolute convergence, we can further use Fubini's theorem to exchange the order of the sums and thus,
\begin{eqnarray*}
\sum_{n=2}^{\infty}w_n\sum_{k=1}^{n-1}\log\left(\frac{1-v_k}{1-v}\right)&=&\sum_{k=1}^{\infty}\log\left(\frac{1-v_k}{1-v}\right)\sum_{n=k+1}^{\infty}w_n\\&=&\sum_{k=1}^{\infty}\log\left(\frac{1-v_k}{1-v}\right)\prod_{n=1}^k(1-v_n)\\
&=&(1-v_1)\log\left(\frac{1-v_1}{1-v}\right)+\sum_{k=2}^{\infty}\log\left(\frac{1-v_k}{1-v}\right)\prod_{n=1}^k(1-v_n).\\
\end{eqnarray*}
Finally, by considering that for all $n$ we define the binary divergence of $v_n$ with respect to $v$ as,
\begin{eqnarray*}
    d(v_n||v)&=&v_n\log\left(\frac{v_n}{v}\right)+(1-v_n)\log\left(\frac{1-v_n}{1-v}\right),
\end{eqnarray*}
the Kullback-Leibler divergence can be finally written as
\begin{eqnarray*}
D_{\mathrm{KL}}(P||P')&=&d(v_1||v)+\sum_{n=2}^{\infty}w_n\log\left(\frac{v_n}{v}\right)+\sum_{n=2}^{\infty}\log\left(\frac{1-v_n}{1-v}\right)\prod_{k=1}^{n}(1-v_k)\nonumber\\
&=&d(v_1||v)+\sum_{n=2}^{\infty}v_n\prod_{k=1}^{n-1}(1-v_k)\log\left(\frac{v_n}{v}\right)+\sum_{n=2}^{\infty}\log\left(\frac{1-v_n}{1-v}\right)\prod_{k=1}^{n}(1-v_k)\nonumber\\
&=&d(v_1||v)+\sum_{n=2}^{\infty}\prod_{k=1}^{n-1}(1-v_k)\left[v_n\log\left(\frac{v_n}{v}\right)+(1-v_n)\log\left(\frac{1-v_n}{1-v}\right)\right]\nonumber\\
&=&d(v_1||v)+\sum_{n=2}^{\infty}\prod_{k=1}^{n-1}(1-v_k)d(v_n||v)\nonumber\\
&=&\sum_{n=1}^{\infty}\prod_{k=1}^{n-1}(1-v_k)d(v_n||v).
\end{eqnarray*}
For the reversed Kullback-Leibler divergence, a similar analysis can be carried out to prove the desired result. In fact, for the geometric case we have directly that $w_n'\leq\lambda^n$ for all $n\geq1$ where $\lambda$ can be chosen as the $\max\{v,1-v\}$.
\end{proof}

\begin{proof}[Proof of Lemma~\ref{expectation_KL}] Taking the expectation of the Kullback-Leibler divergence of the Dirichlet process, $P$, with respect to the geometric process, $P'$ given by equation \eqref{Final_KL} and using Beppo-Levi's theorem yields,
\begin{eqnarray*}
\Esp(D_{\mathrm{KL}}(P||P'))&=&\Esp\left(\sum_{n=1}^{\infty}\left[\prod_{k=1}^{n-1}(1-v_k)\right]\left[v_n\log\left(\frac{v_n}{v}\right)+(1-v_n)\log\left(\frac{1-v_n}{1-v}\right)\right]\right)\\
&=&\sum_{n=1}^{\infty}\Esp\left(\prod_{k=1}^{n-1}(1-v_k)\right)\Esp\left(v_n\log\left(\frac{v_n}{v}\right)+(1-v_n)\log\left(\frac{1-v_n}{1-v}\right)\right)\\
&=&\sum_{n=1}^{\infty}\left(\frac{\theta}{\theta+1}\right)^{n-1}\Esp\left(d(v_n||v)\right)\\
&=&(\theta+1)\Esp\left(d(v_1||v)\right),
\end{eqnarray*}
where
\begin{eqnarray*}
\Esp\left(d(v_1||v)\right)&=&\Esp\left(v_1\log(v_1)\right)-\Esp\left(v_1\log(v)\right)+\Esp\left((1-v_1)\log(1-v_1)\right)-\Esp\left((1-v_1)\log(1-v)\right).
\end{eqnarray*}

Then, to derive a closed expression for this expectation, we use the formulas of the expectation of logarithmically transformed beta random variables provided in Section~\ref{auxiliary_formulas}, from which it can be seen that,
\begin{align}
\Esp(v_1\log(v_1))=&\left(\frac{1}{\theta+1}\right)[\psi(2)-\psi(\theta+2)],\nonumber\\
\Esp(v_1\log(v))=&\left(\frac{1}{\theta+1}\right)[\psi(\alpha)-\psi(\alpha+\beta)],\nonumber\\
\Esp((1-v_1)\log(1-v_1))=&-\frac{\theta}{(\theta+1)^2},\nonumber\\
\Esp((1-v_1)\log(1-v))=&\left(\frac{\theta}{\theta+1}\right)[\psi(\beta)-\psi(\alpha+\beta)]\nonumber.
\end{align}

Hence, by joining these expressions and after some algebraic manipulations we are then able to obtain the desired result, {\it i.e.,},
\begin{eqnarray*}
\Esp(d(v_1||v))&=&\rho[\psi(2)-\psi(\theta+2)-\psi(a)+\psi(a+b)]\nonumber\\&+&(1-\rho)[-\rho-\psi(b)+\psi(a+b)].
\end{eqnarray*}

Now, to obtain the variance we centre our attention on the second moment which is given by
\begin{eqnarray*}
&&\Esp\biggl(\sum_{n=1}^{\infty}\left[\prod_{k=1}^{n-1}(1-v_k)^2\right]d(v_n||v)^2+2\sum_{i<j}\left[\prod_{k=1}^{i-1}(1-v_k)\right]d(v_i||v)\left[\prod_{k=1}^{j-1}(1-v_k)\right]d(v_j||v)\biggr)\\
&=&\Esp\left(\sum_{n=1}^{\infty}\left[\prod_{k=1}^{n-1}(1-v_k)^2\right]d(v_n||v)^2\right)\\
&+&2\Esp\left(\sum_{i<j}\left[\prod_{k=1}^{i-1}(1-v_k)^2\right]\left[\prod_{k=i+1}^{j-1}(1-v_k)\right](1-v_i)d(v_i||v)d(v_j||v)\right)\\
&=&\Esp(d^2(v_1||v))\sum_{n=1}^{\infty}\left(\frac{\theta}{\theta+2}\right)^{n-1}+2\Esp((1-v_1)d(v_1||v)d(v_2||v))\sum_{i<j}\left(\frac{\theta}{\theta+2}\right)^{i-1}\left(\frac{\theta}{\theta+1}\right)^{j-i-1}\\
&=&\left(\frac{\theta+2}{2}\right)\Esp(d^2(v_1||v))+(\theta+2)(\theta+1)\Esp((1-v_1)d(v_1||v)d(v_2||v)),
\end{eqnarray*}
which yields the desired result on the variance. Finally, it only remains to derive a closed expression for the $\Esp(d(v_1||v)^2)$ and for the $\Esp((1-v_1)d(v_1||v)d(v_2||v))$. For the former, we notice after some algebraic manipulations that,
\begin{eqnarray*}
d^2(v_1||v)&=&v_1^2\log^2(v_1)-2v_1^2\log(v_1)\log(v)+v_1^2\log^2(v)\\
&+&(1-v_1)^2\log^2(1-v_1)-2(1-v_1)^2\log(1-v_1)\log(1-v)\\&+&(1-v_1)^2\log^2(1-v)
+2v_1(1-v_1)\log(v_1)\log(1-v_1)\\&-&2v_1(1-v_1)\log(v_1)\log(1-v)
-2v_1(1-v_1)\log(v)\log(1-v_1)\\
&+&2v_1(1-v_1)\log(v)\log(1-v).
\end{eqnarray*}
Hence, it is straightforward to check by using the formulas in Section~\ref{auxiliary_formulas} and the independence of the random variables that
\begin{eqnarray*}
\Esp(d^2(v_1||v))&=&a_1[\Esp(\log^2(X_1))-2\Esp(\log(X_1))\Esp(\log(V))+\Esp(\log^2(V))]\\
&+&a_2\bigl[\Esp(\log^2(1-X_2))-2\Esp(\log(1-X_2))\Esp(\log(1-V))+\Esp(\log^2(1-V))\bigr]\\
&+&a_3[\Esp(\log(X_3)\log(1-X_3))-\Esp(\log(X_3))\Esp(\log(1-V))]\\
&+&a_3[\Esp(\log(V)\log(1-V))-\Esp(\log(V))\Esp(\log(1-X_3))],
\end{eqnarray*}
where $X_1\sim Be(3,\theta)$, $X_2\sim Be(1,\theta+2)$, $X_3
\sim Be(2,\theta+1)$ and the coefficients $a_1,a_2$ and $a_3$ are given by,
\begin{align*}
    a1&=\frac{B(3,\theta)}{B(1,\theta)},&&
    &&a_2=\frac{B(1,\theta+2)}{B(1,\theta)},&&
    &&a_3=\frac{2B(2,\theta+1)}{B(1,\theta)}.
\end{align*}

Now, proceeding in a similar fashion for the $\Esp((1-v_1)d(v_1||v)d(v_2||v))$, observe that
\begin{eqnarray*}
    (1-v_1)d(v_1||v)d(v_2||v)&=&(1-v_1)v_1v_2\log(v_1)\log(v_2)-(1-v_1)v_1v_2\log(v_1)\log(v)\\
    &-&(1-v_1)v_1v_2\log(v_2)\log(v)+(1-v_1)v_1v_2\log^2(v)\\
    &+&(1-v_1)(1-v_2)v_1\log(v_1)\log(1-v_2)\\
    &-&(1-v_1)(1-v_2)v_1\log(v)\log(1-v_j)\\
    &-&(1-v_1)(1-v_2)v_1\log(v_1)\log(1-v)\\
    &+&(1-v_1)(1-v_2)v_1\log(v)\log(1-v)\\
    &+&v_2(1-v_1)^2\log(1-v_1)\log(v_2)-v_2(1-v_1)^2\log(1-v_1)\log(v)\\
    &-&v_2(1-v_1)^2\log(1-v)\log(v_2)+v_2(1-v_1)^2\log(1-v)\log(v)\\
    &+&(1-v_1)^2(1-v_2)\log(1-v_1)\log(1-v_2)\\
    &-&(1-v_1)^2(1-v_2)\log(1-v_1)\log(1-v)\\
    &-&(1-v_1)^2(1-v_2)\log(1-v)\log(1-v_2)+(1-v_1)^2(1-v_2)\log^2(1-v).
\end{eqnarray*}
Hence, by taking expectation we get that the $\Esp((1-v_1)d(v_1||v)d(v_2||v))$ is given by
\begin{eqnarray*}
&&a_1[\Esp(\log(X_1))\Esp(\log(X_2))-\Esp(\log(X_1))\Esp(\log(V))\\
&-&\Esp(\log(X_2))\Esp(\log(V))+\Esp(\log^2(V))]\\
&+&a_2[\Esp(\log(X_1))\Esp(\log(1-X_3))-\Esp(\log(1-X_3))\Esp(\log(V))\\&-&\Esp(\log(X_1))\Esp(\log(1-V))]\\
&+&(a_2+a_3)[\Esp(\log(V)\log(1-V))]\\
&+&a_3[\Esp(\log(1-X_4))\Esp(\log(X_2))-\Esp(\log(1-X_4))\Esp(\log(V))\\&-&\Esp(\log(1-V))\Esp(\log(X_2))]\\
&+&a_4[\Esp(\log(1-X_4))\Esp(\log(1-X_3))-\Esp(\log(1-X_4))\Esp(\log(1-V))]\\
&+&a_4[\Esp(\log(1-X_3))\Esp(\log(1-V))-\Esp(\log^2(1-V))],
\end{eqnarray*}
where $X_1\sim Be(2,\theta+1)$, $X_2\sim Be(2,\theta)$, $X_3\sim Be(1,\theta+1)$, $X_4\sim Be(1,\theta+2)$ and the coefficients $a_1,a_2,a_3$ and $a_4$ are
\begin{align*}
    a1&=\frac{B(2,\theta+1)B(2,\theta)}{B(1,\theta)^2},&&
    &&a_2=\frac{B(2,\theta+1)B(1,\theta+1)}{B(1,\theta)^2},\\
    a_3&=\frac{B(1,\theta+2)B(2,\theta)}{B(1,\theta)^2},&&
    &&a_4=\frac{B(1,\theta+2)B(1,\theta+1)}{B(1,\theta)^2}.
\end{align*}
\end{proof}

\begin{proof}[Proof of Lemma~\ref{expectation_coupled_KL}] Following the same procedure as for the proof of Lemma~\ref{expectation_KL}, we have that the expectation of the Kullback-Leibler divergence for the coupled case, {\it i.e.,}, $v_1=v$ is,
\begin{eqnarray*}
\Esp(D_{\mathrm{KL}}(P||P'))&=&\Esp\left(\sum_{n=1}^{\infty}\left[\prod_{k=1}^{n-1}(1-v_k)\right]\left[v_n\log\left(\frac{v_n}{v_1}\right)+(1-v_n)\log\left(\frac{1-v_n}{1-v_1}\right)\right]\right)\\
&=&\sum_{n=2}^{\infty}\Esp\left(\left[\prod_{k=1}^{n-1}(1-v_k)\right]\left[v_n\log\left(\frac{v_n}{v_1}\right)+(1-v_n)\log\left(\frac{1-v_n}{1-v_1}\right)\right]\right).
\end{eqnarray*}

In this case, we can split the expectation into the following four terms,
\begin{align*}
f_1&=\Esp\left[\prod_{k=1}^{n-1}(1-v_k)\right]\Esp[v_n\log(v_n)],\\
f_2&=\Esp[(1-v_1)\log(v_1)]\Esp\left[\prod_{k=2}^{n-1}(1-v_k)\right]\Esp(v_n),\\
f_3&=\Esp\left[\prod_{k=1}^{n-1}(1-v_k)\right]\Esp[(1-v_n)\log(1-v_n)],\\
f_4&=\Esp[(1-v_1)\log(1-v_1)]\Esp\left[\prod_{k=2}^{n}(1-v_k)\right],
\end{align*}
so that $\Esp(D_{\mathrm{KL}}(P||P'))=f_1-f_2+f_3-f_4$. Then, using the formulas of the expectation of logarithmically transformed beta random variables provided in Section~\ref{auxiliary_formulas}, and after some algebraic manipulations, we obtain that
\begin{align*}
f_1&=\left(\frac{\theta}{\theta+1}\right)^{n-1}\left(\frac{1}{\theta+1}\right)[\psi(2)-\psi(\theta+2)],\\
f_2&=\left(\frac{\theta}{\theta+1}\right)^{n-1}\left(\frac{1}{\theta+1}\right)[\psi(1)-\psi(\theta+2)],\\
f_3&=-\left(\frac{\theta}{\theta+1}\right)^{n-1}\frac{\theta}{(\theta+1)^2},\\
f_4&=-\left(\frac{\theta}{\theta+1}\right)^{n-1}\frac{\theta}{(\theta+1)^2}.
\end{align*}
Then, observe that $f_3=f_4$, so that $\Esp(D_{\mathrm{KL}}(P||P'))=f_1-f_2$. This expression can be further simplified by noting that $\psi(z+1)=\psi(z)+z^{-1}$, therefore,
\begin{eqnarray*}
\Esp(D_{\mathrm{KL}}(P||P'))&=&\sum_{n=2}^{\infty}\left(\frac{\theta}{\theta+1}\right)^{n-1}\left(\frac{1}{\theta+1}\right)[\psi(1)+1-\psi(\theta+2)-\psi(1)+\psi(\theta+2)]\\
&=&\sum_{n=2}^{\infty}\left(\frac{\theta}{\theta+1}\right)^{n-1}\left(\frac{1}{\theta+1}\right)\\
&=&\left(\frac{1}{\theta+1}\right)\sum_{n=2}^{\infty}\left(\frac{\theta}{\theta+1}\right)^{n-1}\\
&=&\frac{\theta}{\theta+1}.
\end{eqnarray*}

Finally, to obtain the variance, it suffices to find the second moment for which we can notice that the first term of the Kullback-Leibler divergence is zero so that the second moment is given by,
\begin{eqnarray*}
&&\Esp\Biggl(\sum_{n=2}^{\infty}\left[\prod_{k=1}^{n-1}(1-v_k)^2\right]d(v_n||v_1)^2
+2\sum_{1<i<j}\left[\prod_{k=1}^{i-1}(1-v_k)\right]d(v_i||v_1)\left[\prod_{k=1}^{j-1}(1-v_k)\right]d(v_j||v_1)\Biggr)\\
&=&\Esp\left(\sum_{n=2}^{\infty}\left[\prod_{k=2}^{n-1}(1-v_k)^2\right](1-v_1)^2d(v_n||v_1)^2\right)\\
&+&2\Esp\left(\sum_{1<i<j}\left[\prod_{k=2}^{i-1}(1-v_k)^2\right]\left[\prod_{k=i+1}^{j-1}(1-v_k)\right](1-v_1)^2(1-v_i)d(v_i||v)d(v_j||v)\right)\\
&=&\Esp((1-v_1)^2d^2(v_2||v_1))\sum_{n=2}^{\infty}\left(\frac{\theta}{\theta+2}\right)^{n-2}\\
&+&2\Esp((1-v_1)^2(1-v_2)d(v_2||v_1)d(v_3||v_2))\sum_{1<i<j}\left(\frac{\theta}{\theta+2}\right)^{i-2}\left(\frac{\theta}{\theta+1}\right)^{j-i-2}\\
&=&\left(\frac{\theta+2}{2}\right)\Esp((1-v_1)^2d^2(v_2||v_1))+(\theta+2)(\theta+1)\Esp((1-v_1)^2(1-v_2)d(v_2||v_1)d(v_3||v_1)),
\end{eqnarray*}
which yields the desired result about the variance. Lastly, we only need to provide closed expressions for the $\Esp((1-v_1)^2d^2(v_2||v_1))$ and for the $\Esp[(1-v_1)^2(1-v_2)d(v_2||v_1)d(v_3||v_1)]$. For the former, we exploit the fact that we already know the value for $d^2(v_2||v_1)$, so that it only remains to multiply it by $(1-v_1)^2$, and take the expectation which yields that,
\begin{eqnarray*}
\Esp((1-v_1)^2d^2(v_2||v_1))&=&a_1[\Esp(\log^2(X_1))-2\Esp(\log(X_1))\Esp(\log(X_2))+\Esp(\log^2(X_2))]\\
&+&2a_2[\Var(\log(1-X_2))]\\
&+&a_3[\Esp(\log(X_3)\log(1-X_3))-\Esp(\log(X_3))\Esp(\log(1-X_2))]\\
&+&a_3[\Esp(\log(X_2)\log(1-X_2))-\Esp(\log(X_2))\Esp(\log(1-X_3))],
\end{eqnarray*}
where $X_1\sim Be(3,\theta)$, $X_2\sim Be(1,\theta+2)$, $X_3
\sim Be(2,\theta+1)$ and the coefficients $a_1,a_2$ and $a_3$ are given by,
\begin{align*}
a1&=\frac{B(3,\theta)B(1,\theta+2)}{B(1,\theta)^2},    &a_2=\left(\frac{B(1,\theta+2)}{B(1,\theta)}\right)^2,
&a_3=\frac{2B(2,\theta+1)B(1,\theta+2)}{B(1,\theta)^2}.
\end{align*}

Finally, for the $\Esp((1-v_1)^2(1-v_2)d(v_2||v_1)d(v_3||v_1))$, we already know an expression for $(1-v_2)d(v_2||v_1)d(v_3||v_1)$, so that it only remains to multiply by $(1-v_1)^2$. Then, by taking the expectation, it is possible to show, after some algebraic manipulations that,
\begin{eqnarray*}
\Esp((1-v_1)^2(1-v_2)d(v_2||v_1)d(v_3||v_1))
&=&a_1[\Esp(\log(X_1))\Esp(\log(X_2))-\Esp(\log(X_3))\Esp(\log(X_1))\\&-&\Esp(\log(X_3))\Esp(\log(X_2))+\Esp(\log^2(X_3))]\\
&+&a_2[\Esp(\log(X_1))\Esp(\log(1-X_4))\\
&-&\Esp(\log(1-X_4))\Esp(\log(X_3))\\
&-&\Esp(\log(X_1))\Esp(\log(1-X_3))]\\
&+&(a_2+a_3)[\Esp(\log(X_3)\log(1-X_3))]\\
&-&a_3[\Esp(\log(X_3))\Esp(\log(1-X_3))]\\&+&a_4\Var(\log(1-X_3)),
\end{eqnarray*}
where $X_1\sim Be(2,\theta+1)$, $X_2\sim Be(2,\theta)$, $X_3\sim Be(1,\theta+2)$, $X_4\sim Be(1,\theta+1)$ and the coefficients $a_1,a_2,a_3$ and $a_4$ are given by,
\begin{align*}
a1&=\frac{B(1,\theta+2)B(2,\theta+1)B(2,\theta)}{B(1,\theta)^3},&&    a_2=\frac{B(1,\theta+2)B(2,\theta+1)B(1,\theta+1)}{B(1,\theta)^3},\\
a_3&=\frac{B(1,\theta+2)^2B(2,\theta)}{B(1,\theta)^3},&&    a_4=\frac{B(1,\theta+2)^2B(1,\theta+1)}{B(1,\theta)^3}.
\end{align*}
\end{proof}

\subsection{Proofs of section 4}
\label{proofsSec4}

\begin{proof}[Proof of Theorem~\ref{KL_SBP_GEOM}]
The proof of this theorem follows the same reasoning as Theorem 1. For this, we first appreciate that if we let $P_{\beta}$ be the exchangeable stick-breaking process, with length variables driven by a Dirichlet process $(\beta,\nu_0)$ with base measure $\nu_0=Be(1,\theta)$, and $P'$ the geometric process $(1,\theta)$, then the Kullback-Leibler divergence of $P_{\beta}$ with respect to $P'$ is the same as \eqref{DKL}, that is,
\begin{eqnarray*}
    D_{\mathrm{KL}}(P_{\beta}||P')&=&v_1\log\left(\frac{v_1}{v}\right)+\sum_{n=2}^{\infty}\log\left(\frac{v_n}{v}\right)+\sum_{n=2}^{\infty}w_n\sum_{k=1}^{n-1}\log\left(\frac{1-v_k}{1-v}\right).
\end{eqnarray*}
Then, we need to show the almost surely absolute convergence of the series. For this, a routine application of de Finetti's representation theorem allows one to see, that marginally the length-variables follow a $Be(1,\theta)$ distribution, so that \eqref{eq:AlmostSureConvergence2} holds. Furthermore, if we let $\mu=\Esp(\log(1-V))$ where $V\sim Be(1,\theta)$, then 
\begin{eqnarray*}
    \Pro\left(\lim_{n\to\infty}\frac{1}{n}\sum_{k=1}^n\log(1-v_k)=\mu\right)=\int\Pro\left(\lim_{n\to\infty}\frac{1}{n}\sum_{k=1}^n\log(1-v_k)=\mu\Big|\nu\right)Q(d\nu),
\end{eqnarray*}
where $\Pro\left(\lim_{n\to\infty}\frac{1}{n}\sum_{k=1}^n\log(1-v_k)=\mu\Big|\nu\right)=1$ since, conditioned on $\nu$, the length-variables are iid so that we can use the strong law of the large numbers as we did in \eqref{eq:AlmostSureConvergence1}. Hence, 
\begin{eqnarray*}
    w_n\leq\left(e^{\Esp[\log(1-V)+\epsilon]}\right)^n=\lambda^n,
\end{eqnarray*}
and since $\nu_0$ places zero mass to $\{0\}$ we have a proper species-sampling model, that is, $\sum_{n=1}^{\infty}w_n=1$ a.s. Therefore, we have the necessary conditions to ensure the almost surely convergence of the series and it that way we obtain that,
\begin{eqnarray*}
D_{\mathrm{KL}}(P_{\beta}||P')&=&\sum_{n=1}^{\infty}\left[\prod_{k=1}^{n-1}(1-v_k)\right]d(v_n||v),
\end{eqnarray*}
as desired.
\end{proof}

\begin{proof}[Proof of Lemma~\ref{exp_exch_sbp}]
Consider the Kullback-Leibler divergence of $P_{\beta}$ with respect to $P$, which is given by,
\begin{eqnarray*}
    D_{\mathrm{KL}}(P_{\beta}||P)&=&\sum_{n=1}^{\infty}\prod_{j=1}^{n-1}(1-v_j)d(v_n||v).
\end{eqnarray*}
By Beppo-Levi's theorem, we have that
\begin{eqnarray*}
    \Esp(D_{\mathrm{KL}}(P_{\beta}||P'))&=&\sum_{n=1}^{\infty}\Esp\left(\prod_{j=1}^{n-1}(1-v_j)d(v_n||v)\right).
\end{eqnarray*}
Now, if for all $n\geq1$ we define, 
\begin{eqnarray*}
    f(v_1,\ldots,v_n)&=&\prod_{j=1}^{n-1}(1-v_j)d(v_n||v),
\end{eqnarray*}
we can then obtain its expectation using the formula described in \cite{Leyva_2021}, that establishes that for a measurable function $f$ of exchangeable random variables, its expected value is given by,
\begin{eqnarray*}
    \Esp(f(v_1,...,v_n))&=&\sum_{\{A_1,...,A_k\}}\left[\int f(v_{l_1},...,v_{l_n})\prod_{j=1}^k\prod_{r\in A_j}\Ind_{(l_r=j)}\nu_0(dv_1)\cdots\nu_0(dv_k)\right]p(n_1,\ldots,n_k),
\end{eqnarray*}
where the sum ranges over all the partitions $\{A_1,\ldots,A_k\}$, and $p(\cdot)$ is the exchangeable partition probability function (EPPF) associated to the nonparametric prior chosen. For ease of notation, if we denote for all $n\geq1$ the set of all partitions of $\{1,\ldots,n\}$ as $\mathscr{P}([n])$, and $\{v^*_1,\ldots,v^*_{k}\}$ the $k$ distinct values found in the set $\{v_1,\ldots,v_n\}$, then the expectation of Kullback-Leibler divergence can be written as
\begin{eqnarray*}
    \Esp(D_{\mathrm{KL}}(P_{\beta}||P'))&=&\sum_{n=1}^{\infty}\sum_{\pi\in\mathscr{P}([n])}F_{\theta}(\pi)p_{\beta}(\pi),
\end{eqnarray*}
where
\begin{eqnarray*}
F_{\theta}(\pi)=\Esp\left(\left(\prod_{j=1}^{k-1}(1-v_j)^{n_j}\right)(1-v_{k})^{n_k-1}d(v_k||v)\right).
\end{eqnarray*}
Hence, if $p_{\beta}(\cdot)$ is the EPPF of the Dirichlet process $(\beta,\nu_0)$, we have that,
\begin{eqnarray*}
 &&F_{\theta}(\pi)p_{\beta}(\pi)\\&=&\Esp\left(\left(\prod_{j=1}^{k-1}(1-v_j)^{n_j}\right)(1-v_{k_n})^{n_k-1}d(v_k||v)\right)\frac{\beta^{k}}{(\beta)^{(n)}}\prod_{j=1}^{k}(n_j-1)!,
\end{eqnarray*}
as required.
\end{proof}

\begin{proof}[Proof of Theorem~\ref{bound_cont_comp}]
First, we centre our attention on the upper bound, so that without loss of generality, in the sequel we assume that any  $\pi = \{A_{1},\ldots,A_{k}\} \in\mathcal{P}([n])$ satisfies that $1\in A_{1}$. Furthermore, we also assume that $n\in A_{k}$ whenever $n \not\in A_{1}$. Note that if $n\in A_{1}$, then $F_{\theta}(\pi) = 0$. Also, observe that if $n\not\in A_{1}$ then,
\begin{equation*}
    F_{\theta}(\pi) = \prod_{j=2}^{k-1} \frac{\theta}{\theta + n_j} \times \mathbb{E}((1-V)^{n_1} (1-V')^{ n_k  - 1} d(V' \Vert V)).
\end{equation*}
The inequality $(1-V)^{n_1} (1-V')^{n_k - 1} \leq (1-V)$, and the fact that $\mathbb{E}((1-V)d(V' \Vert V)) = \theta/(\theta+1)^{2}$ lead to
\begin{equation*}
\label{aux_cont}
    F_{\theta}(\pi) \leq \frac{1}{\theta} \bigg(\frac{\theta}{\theta+1}\bigg)^{k}.
\end{equation*}
For ease of notation, we let $\alpha = \theta/(\theta+1)$. Then the previous inequality implies that
\begin{align}
\label{eq:ProofTheoremComplexityPerformanceTradeoffMainBound}
     F_{\theta}(\pi) p_{\beta}(\pi) \leq \frac{1}{\theta} \frac{(\alpha\beta)^{k}}{\beta^{(n)}} \prod_{j=1}^{k} (n_j - 1)!
     = \frac{1}{\theta} \frac{(\alpha\beta)^{(n)}}{\beta^{(n)}} \frac{(\alpha\beta)^{k}}{(\alpha\beta)^{(n)}} \prod_{j=1}^{k} (n_j - 1)!.
\end{align}
Since $F_{\theta}(1_{n}) = 0$, we have that
\begin{equation*}
    \sum_{\pi\in\mathcal{P}([n])} F_{\theta}(\pi) p_{\beta}(\pi) = \sum_{\pi<1_{n}} F_{\theta}(\pi) p_{\beta}(\pi).
\end{equation*}
The previous equation and 
\eqref{eq:ProofTheoremComplexityPerformanceTradeoffMainBound} imply that,
\begin{align*}
    \sum_{\pi\in\mathcal{P}([n])} F_{\theta}(\pi) p_{\beta}(\pi) &\leq \frac{1}{\theta} \frac{(\alpha\beta)^{(n)}}{\beta^{(n)}} \sum_{\pi<1_{n}} \frac{(\alpha\beta)^{k}}{(\alpha\beta)^{(n)}} \prod_{j=1}^{k} (n_j - 1)!\\
    &= \frac{1}{\theta} \frac{(\alpha\beta)^{(n)}}{\beta^{(n)}} \sum_{\pi\in\mathcal{P}([n])} \frac{(\alpha\beta)^{k}}{(\alpha\beta)^{(n)}} \prod_{j=1}^{k} (n_j - 1)! - \frac{\beta}{\theta+1} \frac{(n-1)!}{\beta^{(n)}}\\
    &= \frac{1}{\theta} \frac{(\alpha\beta)^{(n)}}{\beta^{(n)}} - \frac{\beta}{\theta+1} \frac{(n-1)!}{\beta^{(n)}}.
\end{align*}
In particular, we have that
\begin{equation*}
    D_{\theta}(\beta) \leq \frac{1}{\theta} \sum_{n=1}^{\infty} \frac{(\alpha\beta)^{(n)}}{\beta^{(n)}} - \frac{\beta}{\theta+1} \sum_{n=1}^{\infty} \frac{(n-1)!}{\beta^{(n)}}.
\end{equation*}
Therefore, Lemma~\ref{Lemma:SumRisingFactorials} and Lemma~\ref{Lemma:SumQuotientRisingFactorials} lead to
\begin{align*}
    D_{\theta}(\beta) \leq \frac{\beta}{\theta+1} \bigg(\frac{1}{\beta/(\theta+1)-1} - \frac{1}{\beta-1}\bigg)= \frac{\theta}{\theta+1} \frac{\beta^{2}}{(\beta-1)(\beta-(\theta+1))},
\end{align*}
as required. Then to prove that the expectation is continuous, first let $\lambda>1$ and, using the same reasoning as before, we can see that
\begin{eqnarray*}
    D_{\theta}(\lambda\beta)&\leq&\sum_{n=1}^{\infty}\sum_{\pi\in\mathcal{P}([n])}\frac{1}{\theta} \frac{(\alpha\lambda\beta)^{k}}{(\lambda\beta)^{(n)}} \prod_{j=1}^{k} (n_j - 1)!\\
&\leq&\sum_{n=1}^{\infty}\sum_{\pi\in\mathcal{P}([n])}\frac{1}{\theta} \frac{(\alpha\lambda\beta)^{k}}{(\beta)^{(n)}} \prod_{j=1}^{k} (n_j - 1)!\\
    &=&\frac{1}{\theta}\sum_{n=1}^{\infty}\frac{(\alpha\lambda\beta)^{(n)}}{(\beta)^{(n)}}\sum_{\pi\in\mathcal{P}([n))} \frac{(\alpha\lambda\beta)^{k}}{(\alpha\lambda\beta)^{(n)}} \prod_{j=1}^{k} (n_j - 1)!\\
    &=&\frac{1}{\theta}\sum_{n=1}^{\infty}\frac{(\alpha\lambda\beta)^{(n)}}{(\beta)^{(n)}}.
\end{eqnarray*}
Then, from Lemma~\ref{Lemma:SumQuotientRisingFactorials}, $D_{\theta}(\lambda\beta)$ is finite if and only if $(1-\lambda\alpha)\beta>1$ and $(1-\alpha\lambda)\in(0,1)$, which yields that, for fixed $\theta$ and $\beta>\theta+1$, $\lambda$ must satisfy that $1<\lambda<\frac{\theta+1}{\theta}\left(1-\frac{1}{\beta}\right)=\lambda_0$. Now, observe that if we let
\begin{eqnarray*}
    f_n(\lambda)&=&\sum_{\pi\in\mathcal{P}([n])}F_{\theta}(\pi)p_{\lambda\beta}(\pi),
\end{eqnarray*}
and take $\lambda_0^*=\frac{1+\lambda_0}{2}$, then for all $\lambda\in(1,\lambda_0^*)$,
\begin{align*}
    |f_n(\lambda)|\leq g_n(\lambda_0^*)=\sum_{\pi\in\mathcal{P}([n])}\frac{1}{\theta} \frac{(\alpha\lambda_0^*\beta)^{k}}{(\beta)^{(n)}} \prod_{j=1}^{k} (\lvert A_{j} \rvert - 1)!.
\end{align*}
Note that $g_n(\lambda_0^*)$ is an integrable function with respect to the counting measure. Hence, an immediate application of the dominated convergence theorem allows us to guarantee that
\begin{eqnarray*}
    \lim_{\lambda\to1^+}D_{\theta}(\lambda\beta)&=&D_{\theta}(\beta).
\end{eqnarray*}
Now, in a similar manner if we take $\lambda<1$ we have 
\begin{eqnarray*}
    D_{\theta}(\lambda\beta)&\leq&\sum_{n=1}^{\infty}\sum_{\pi\in\mathcal{P}([n])}\frac{1}{\theta} \frac{(\alpha\beta)^{k}}{(\lambda\beta)^{(n)}} \prod_{j=1}^{k} (n_j - 1)!\\
    &=&\frac{1}{\theta}\sum_{n=1}^{\infty}\frac{(\alpha\beta)^{(n)}}{(\lambda\beta)^{(n)}}\sum_{\pi\in\mathcal{P}([n])} \frac{(\alpha\beta)^{k}}{(\alpha\beta)^{(n)}} \prod_{j=1}^{k} (n_j - 1)!\\
    &=&\frac{1}{\theta}\sum_{n=1}^{\infty}\frac{\left(\frac{\alpha}{\lambda}\lambda\beta\right)^{(n)}}{(\lambda\beta)^{(n)}},
\end{eqnarray*}
which is a convergent series whenever $\alpha<\lambda$. Then, in this case we let $\lambda_0^*=\frac{1+\alpha}{2}$ so that, for all $\lambda\in(\lambda_0^*,1)$,
\begin{align*}
    |f_n(\lambda)|\leq g_n(\lambda_0^*)=\sum_{\pi\in\mathcal{P}([n])}\frac{1}{\theta} \frac{(\alpha\beta)^{k}}{(\lambda_0^*\beta)^{(n)}} \prod_{j=1}^{k} (n_j - 1)!.
\end{align*}
Note that in this case, $g_n(\lambda_0^*)$ is also an integrable function with respect to the counting measure. Hence, an immediate application of the dominated convergence theorem allows us to guarantee that
\begin{eqnarray*}
    \lim_{\lambda\to1^-}D_{\theta}(\lambda\beta)&=&D_{\theta}(\beta).
\end{eqnarray*}
Therefore, we get the desired result on the continuity of $D_{\theta}(\beta)$. Finally, it is immediate to notice that if $p(\cdot)$ is the EPPF of the Dirichlet process we have that for all $n\geq 1$,
\begin{align*}
    \lim_{\beta\to\infty} \left[\frac{\beta^{k}}{(\beta)^{(n)}}\prod_{j=1}^{k}(n_j-1)!\right]=
    \begin{cases}
    1&\text{if } k=n,\\
    0&\text{if } k<n.
    \end{cases}
\end{align*}
Therefore, for all $n\geq1$ and $\pi\in\mathscr{P}([n])$, we obtain that
\begin{eqnarray*}
    \lim_{\beta\to\infty}p_{\beta}(\pi)&=&\delta_{(0_n)}(\pi),
\end{eqnarray*}
where $\delta$ is the Dirac measure, and $0_n\in\mathscr{P}([n])$ is the partition induced when there are no ties among the length variables. Hence, due to the continuity and the fact that $f_n$ is dominated by $g_n$ we get that,
\begin{eqnarray*}
    \lim_{\beta\to\infty}D_{\theta}(\beta)    &=&\sum_{n=1}^{\infty}\sum_{\sigma\in\mathscr{P}([n])}\lim_{\beta\to\infty}F_{\theta}(\sigma)\pi_{\beta}(\sigma)\\
    &=&\sum_{n=1}^{\infty}F_{\theta}(0_n)\\&=&\sum_{n=1}^{\infty}\Esp\left(\prod_{j=1}^{n-1}(1-v_j)d(v_n||v)\right),
\end{eqnarray*}
where the last term is $\Esp(D_{\mathrm{KL}}(P||P'))$ as stated.

\end{proof}

\begin{proof}[Proof of Lemma~\ref{Lemma:SumRisingFactorials}]
The convergence of the series for $\beta>1$ follows immediately from Gautschi's inequality. Indeed, if we consider $\beta=1+p$, with $p\in(0,1)$, then
\begin{align*}
    \frac{(n-1)!}{\beta^{(n)}}=\frac{\varGamma(1+p)\varGamma(n+1)}{n(n+p)\varGamma(n+p)}<\frac{\varGamma(1+p)(n+1)^{1-p}}{n(n+p)}\sim \frac{1}{n^{1+p}},
\end{align*}
which is convergent if and only if $p>0$, and hence $\beta>1$. To obtain a closed expression for this series, we use the following series expansion of $\beta^{-(m+1)}$ for $m$ a positive integer and $Re(\beta)>0$ \cite[\textit{see e.g.}][]{tweddle2012james}
\begin{eqnarray*}
    \frac{1}{\beta^{m+1}}&=&\sum_{n=0}^{\infty}\frac{\sigma(n+m,m)}{\beta(\beta+1)\cdots(\beta+n+m)},
\end{eqnarray*}
where $\sigma(n+m,m)$ is the unsigned Stirling number of the first kind. Hence, if we choose $m=1$ we have that,
\begin{eqnarray*}
    \frac{1}{\beta}&=&\sum_{n=0}^{\infty}\frac{\sigma(n+1,1)}{(\beta+1)\cdots(\beta+n+1)},
\end{eqnarray*}
where $\sigma(n+1,1)=n!$. Therefore, for $\beta-1>0$ we have that
\begin{eqnarray*}
    \frac{1}{\beta-1}&=&\sum_{n=0}^{\infty}\frac{n!}{(\beta)\cdots(\beta+n)}=\sum_{n=1}^{\infty}\frac{(n-1)!}{\beta(\beta+1)\cdots(\beta+n-1)},
\end{eqnarray*}
as required.
\end{proof}

\begin{proof}[Proof of Lemma~\ref{Lemma:SumQuotientRisingFactorials}]
For this series, we use the following result found in \cite{norlund}, that establishes that for $Re(z)>w$ we have the following expansion series for $(z-w)^{-1}$
\begin{eqnarray*}
\frac{1}{z-w}&=&\sum_{n=0}^{\infty}\frac{w(w+1)\cdots(w+n-1)}{z(z+1)\cdots(z+n)}.
\end{eqnarray*}
Now, we notice that,
\begin{eqnarray*}
    \sum_{n=1}^{\infty}\frac{(\lambda\beta)^{(n)}}{(\beta)^{(n)}}&=&\sum_{n=1}^{\infty}\frac{(\lambda\beta)(\lambda\beta+1)(\lambda\beta+2)\cdots(\lambda\beta+n-1)}{(\beta)(\beta+1)(\beta+2)\cdots(\beta+n-1)}\\
    &=&\lambda\beta\sum_{n=1}^{\infty}\frac{(\lambda\beta+1)([\lambda\beta+1]+1)\cdots([\lambda\beta+1]+n-2)}{(\beta)(\beta+1)(\beta+2)\cdots(\beta+n-1)}\\
    &=&\lambda\beta\sum_{n=0}^{\infty}\frac{(\lambda\beta+1)([\lambda\beta+1]+1)\cdots([\lambda\beta+1]+n-1)}{(\beta)(\beta+1)(\beta+2)\cdots(\beta+n)}\\
    &=&\frac{\lambda\beta}{\beta-\lambda\beta-1},
\end{eqnarray*}
where the convergence is achieved whenever $\beta>\lambda\beta+1$ as required.
\end{proof}
\end{appendix}

\clearpage

\end{document}